\documentclass{amsart}

\usepackage{amsmath,amssymb,amsthm,mathtools}
\usepackage{aliascnt}
\usepackage{amsrefs}
\usepackage{enumitem}
\usepackage{microtype}
\usepackage[dvipsnames]{xcolor}
\usepackage[colorlinks=true,linkcolor=blue!55!black,citecolor=green!45!black,urlcolor=red!60!black]{hyperref}
\numberwithin{equation}{section}

\usepackage{fourier}
\usepackage[margin=1.3in]{geometry}

\newtheorem{theorem}{Theorem}[section]
\newaliascnt{proposition}{theorem}
\newtheorem{proposition}[proposition]{Proposition}
\aliascntresetthe{proposition}
\newaliascnt{lemma}{theorem}
\newtheorem{lemma}[lemma]{Lemma}
\aliascntresetthe{lemma}
\newaliascnt{corollary}{theorem}
\newtheorem{corollary}[corollary]{Corollary}
\aliascntresetthe{corollary}
\newaliascnt{assumption}{theorem}

\aliascntresetthe{assumption}
\theoremstyle{definition}
\newaliascnt{definition}{theorem}

\aliascntresetthe{definition}
\newaliascnt{remark}{theorem}
\newtheorem{remark}[remark]{Remark}
\aliascntresetthe{remark}
\newaliascnt{example}{theorem}

\aliascntresetthe{example}
\usepackage[capitalise,noabbrev]{cleveref}
\crefname{theorem}{Theorem}{Theorems}
\Crefname{theorem}{Theorem}{Theorems}
\crefname{proposition}{Proposition}{Propositions}
\Crefname{proposition}{Proposition}{Propositions}
\crefname{lemma}{Lemma}{Lemmas}
\Crefname{lemma}{Lemma}{Lemmas}
\crefname{corollary}{Corollary}{Corollaries}
\Crefname{corollary}{Corollary}{Corollaries}
\crefname{assumption}{Assumption}{Assumptions}
\Crefname{assumption}{Assumption}{Assumptions}
\crefname{definition}{Definition}{Definitions}
\Crefname{definition}{Definition}{Definitions}
\crefname{remark}{Remark}{Remarks}
\Crefname{remark}{Remark}{Remarks}
\crefname{example}{Example}{Examples}
\Crefname{example}{Example}{Examples}
\crefname{appendix}{Appendix}{Appendices}
\Crefname{appendix}{Appendix}{Appendices}

\DeclareMathOperator{\Ent}{Ent}

\newcommand{\ud}{\,\mathrm{d}}
\newcommand{\RR}{\mathbb{R}}
\newcommand{\e}{\mathrm{e}}

\newcommand{\La}{\mathcal L_a}
\newcommand{\Ls}{\mathcal L_s}
\newcommand{\Pt}{\mathcal P}
\newcommand{\Piv}{\Pi_v}
\newcommand{\Entv}{\Ent_v}
\newcommand{\Iv}{I_v}
\newcommand{\muT}{\bar\mu_T}
\newcommand{\mux}{\mu_x}
\newcommand{\kv}{{\varphi}}

\DeclareMathOperator{\divg}{div}

\DeclareMathOperator{\tr}{tr}
\newcommand{\PPtwo}{\mathcal P_2}
\newcommand{\Cs}{\mathcal C_\sigma}
\newcommand{\Ds}{\mathcal D_\sigma}
\newcommand{\Id}{\mathrm{Id}}

\usepackage{todonotes}

\title[Space-Time HWI and the Wasserstein Current]{Space-Time HWI Inequalities and the Wasserstein Current for Kinetic Equations}
\thanks{B.L. is supported by the National Key R\&D Program of China grant 2024YFA1016000, the City University of Hong Kong start-up fund 7200843, Hong Kong
RGC Early Career Scheme grant 21301326, and NSFC grant 12601817. The work of J.L. is supported in part by the National Science Foundation under grant DMS-2309378. J.L. thanks Felix Otto for some insightful suggestions.}
\author{Bowen Li}
\address{Department of Mathematics, City University of Hong Kong}
\email{boweli4@cityu.edu.hk}

\author{Jianfeng Lu}
\address{Department of Mathematics, Department of Physics, and Department of Chemistry, Duke University}
\email{jianfeng@math.duke.edu}
\date{}
\hypersetup{pdftitle={Space-Time HWI Inequalities and the Wasserstein Current for Kinetic Equations},pdfauthor={Bowen Li and Jianfeng Lu},pdfsubject={Kinetic transport, relative entropy, and space-time functional inequalities}}

\begin{document}
\begin{abstract}
We prove a space-time analogue of the HWI inequality of Otto-Villani for controlled kinetic equations with a semiconvex confining potential. It bounds the time-averaged relative entropy by spatial and conditional velocity transport distances, the velocity Fisher information, an intrinsic action of the path, and a boundary term given by the time endpoint values of the Wasserstein current. The key step is a differential comparison inequality for the Wasserstein current, which captures the transfer of dissipation from velocity to position. Along the Langevin semigroup, this inequality yields a kinetic evolution variational inequality involving the second time derivative of the squared spatial Wasserstein distance and a damping term. Under a compatible spatial Talagrand's transport inequality, we obtain a boundary-free space-time logarithmic Sobolev inequality on every interval $[0,T]$, $T>0$, with the optimal $T^{-2}$ short-time scaling of the constant. As a consequence, we also establish entropy decay for every friction coefficient and every initial density of finite entropy, with optimal cubic order at short times and exponential decay at long times. 
\end{abstract}

\maketitle

\section{Introduction}\label{sec:introduction}

The HWI inequality of Otto and Villani~\cite{OttoVillani2000} connects relative entropy, the quadratic Wasserstein distance, and Fisher information. Specifically, if $\mux(\ud x)=Z_x^{-1}\e^{-U(x)}\ud x$ and $\nabla^2U\succeq K\Id$ for some $K \in \mathbb{R}$, then a probability density $q$ relative to $\mux$ satisfies
\begin{equation}\label{eq:classical-hwi}
 \Ent_{\mux}(q)
 \le W_2(q\mux,\mux)\left(\int\frac{|\nabla_xq|^2}{q}\ud\mux\right)^{1/2}
       -\frac K2W_2^2(q\mux,\mux).
\end{equation}
Here $\Ent_{\mux}(q)=\int q\log q\ud\mux$ is the entropy, and $W_2$ denotes the quadratic Wasserstein distance. If $K>0$, completing the square shows that the right-hand side of \eqref{eq:classical-hwi} is bounded above by \(\frac{1}{2K}\int \frac{|\nabla q|^2}{q}\,\ud\mu_x\), yielding the logarithmic Sobolev inequality (LSI):
\begin{equation}\label{eq:classical-lsi}
 \Ent_{\mux}(q)\le\frac1{2K}\int\frac{|\nabla_xq|^2}{q}\ud\mux.
\end{equation}
If $K\le0$ and $\mux$ satisfies the  Talagrand's transport inequality $W_2^2(q\mux,\mux)\le\frac2\rho\Ent_{\mux}(q)$ with $\rho>0$ and $\rho+K>0$, then \eqref{eq:classical-hwi} gives
\begin{equation}\label{eq:classical-lsi-negative}
 \Ent_{\mux}(q)\le\frac{2\rho}{(\rho+K)^2}\int\frac{|\nabla_xq|^2}{q}\ud\mux.
\end{equation}
The HWI and LSI fit into the variational interpretation of diffusion equations developed in~\cites{JordanKinderlehrerOtto1998,Otto2001}, known as the Otto calculus. In this perspective, the overdamped equation $\partial_tq=\Delta_xq-\nabla U\cdot\nabla_xq$ is the Wasserstein gradient flow of relative entropy. For a fixed probability measure $\sigma$ with finite second moment and finite $\Ent(\sigma\mid\mux)$, the evolution variational inequality (EVI) reads~\cites{OttoVillani2000,ambrosio2005gradient}
\begin{equation}\label{eq:overdamped-evi-intro}
\frac12\frac{\ud}{\ud t}W_2^2(q_t\mu_x,\sigma)
\le
\Ent(\sigma\mid\mu_x)-\Ent_{\mu_x}(q_t)
-\frac K2W_2^2(q_t\mu_x,\sigma).
\end{equation}
The HWI inequality \eqref{eq:classical-hwi}, and hence the LSI, can be recovered from \eqref{eq:overdamped-evi-intro} with $\sigma=\mux$: the Wasserstein first-variation formula expresses the negative time derivative of $\frac12W_2^2(q_t\mux,\mux)$ as a pairing of the entropy gradient with the optimal displacement, which is then bounded by Cauchy--Schwarz; see~\cite{ambrosio2005gradient}. 

This work concerns the underdamped Langevin dynamics. With friction $\gamma>0$, the infinitesimal generator is given by $\La+\gamma\Ls$, where the Hamiltonian transport and velocity Ornstein--Uhlenbeck generators are
\begin{equation}\label{eq:generators-intro}
 \La=v\cdot\nabla_x-\nabla U\cdot\nabla_v,
 \qquad \Ls=\Delta_v-v\cdot\nabla_v.
\end{equation}
The Hamiltonian operator couples position and velocity, while $\gamma\Ls$ accounts for friction and diffusion in velocity. These dynamics preserve the Gibbs measure
\begin{equation}\label{eq:gibbs-intro}
 \mu(\ud x\ud v)=\mux(\ud x)\,\kv(\ud v),
 \qquad
 \kv(\ud v)=(2\pi)^{-d/2}\e^{-|v|^2/2}\ud v.
\end{equation}
For densities $g_t$ relative to $\mu$, the evolution equation and its semigroup are
\begin{equation} \label{eq:langevinpath}
     (\partial_t+\La)g_t=\gamma\Ls g_t,
 \qquad g_t=\Pt_tg_0,
 \qquad \Pt_t=\e^{t(-\La+\gamma\Ls)}.
\end{equation}
The entropy dissipation identity is
\begin{equation}\label{eq:entropy-dissipation-intro}
 -\frac{\ud}{\ud t}\Ent_\mu(\Pt_tg)
       =\gamma \Iv(\Pt_tg),
 \qquad \Iv(g)=\int\frac{|\nabla_vg|^2}{g}\ud\mu.
\end{equation}
An instantaneous bound of entropy by $\Iv$ is impossible: a nonconstant density depending only on $x$ has zero velocity Fisher information. The missing coercivity must therefore be supplied by Hamiltonian transport over a time interval. The transfer of coercivity from velocity to position is the basic difficulty addressed by hypocoercivity; see~\cite{Villani2009}.

We define the \emph{Wasserstein current} relative to the spatial equilibrium $\mu_x$, from the first-variation formula of $W^2_2(q \mu_x,\mu_x)$, by
\begin{equation}\label{eq:cot-intro}
 \mathcal C_{\mux}(g)=\int j(x)\cdot\bigl(x-T_q(x)\bigr)\ud\mux(x),
 \qquad q=\Piv g,\quad j=\Piv(vg),
\end{equation}
where $\Piv$ denotes integration against $\kv$ and $T_q$ is the optimal transportation map from $q\mux$ to $\mux$. A negative Wasserstein current means that the spatial law is moving toward equilibrium in $W_2$; a positive current means it is moving away. In our earlier entropy-decay argument~\cite{Lu2026}, the Wasserstein current, denoted there by $\mathcal C_{\rm OT}$, serves as a corrector that turns the entropy into a Lyapunov function.
In~\cite{LiLu2026}, we subsequently proved a space-time LSI for controlled kinetic paths and used an interpolation between the forward and backward semigroups to establish hypocoercive hypercontractivity. We follow the notation of these works: a measurable velocity current $W_t(x,v)\in\RR^d$ drives a density path $g$ if
\begin{equation}\label{eq:controlled-intro}
       (\partial_t+\La)g=\nabla_v^*W,
       \qquad \nabla_v^*W=-\divg_vW+v\cdot W.
\end{equation}
The intrinsic action of a path measures its departure from Hamiltonian transport, averaged over $[0,T]$:
\begin{equation}\label{eq:action-intro}
 \mathcal A_T(g):=\inf_W\frac1T\int_0^T\!\int\frac{|W_t|^2}{g_t}\ud\mu\ud t,
\end{equation}
where the infimum is over all currents satisfying \eqref{eq:controlled-intro}.

\smallskip

The space-time LSI of~\cite{LiLu2026} is the kinetic counterpart of \eqref{eq:classical-lsi}. In the overdamped setting, the LSI follows from the HWI inequality \eqref{eq:classical-hwi} and the EVI \eqref{eq:overdamped-evi-intro} behind it. The goal of this work is to identify the kinetic counterparts of these two inequalities. We prove a space-time HWI inequality for controlled paths, from which space-time LSIs follow. We also prove a kinetic evolution variational inequality along the Langevin semigroup, which involves the second time derivative of the squared spatial Wasserstein distance and yields entropy decay.

\subsection*{Main results}
We now state our main results. Let $\PPtwo$ be the space of probability measures with finite second moment, equipped with $W_2$.  We assume that the potential is smooth and has a uniform lower Hessian bound, and that the associated Gibbs measure \eqref{eq:gibbs-intro} has finite second moment. More precisely, for some $K\in\RR$,
\begin{equation}\label{eq:intro-potential}
 U\in C^\infty(\RR^d),\qquad
 \nabla^2U\succeq K\Id, \qquad
 \mu\in\PPtwo(\RR^{2d}),
\end{equation}
with $Z_x = \int_{\RR^d}\e^{-U(x)}\ud x < \infty$. Note that no upper bound on $\nabla^2U$ is imposed. 

Let $g_t$ be probability densities with respect to $\mu$ such that $t\mapsto g_t\mu$ is narrowly continuous on $[0,T]$. A measure-valued path is narrowly continuous if its integral against every bounded continuous test function is continuous in time. In addition we assume $g_0\mu\in\PPtwo(\RR^{2d})$, and suppose that \eqref{eq:controlled-intro} holds in distributions with
\begin{equation}\label{eq:intro-path}
 \Ent_\mu(g_0)<\infty,\qquad
 \int_0^T\!\int\frac{|W_t|^2}{g_t}\ud\mu\ud t<\infty,
 \qquad \int_0^T \Iv(g_t)\ud t<\infty.
\end{equation}
The quotient $|W_t|^2/g_t$ is zero when the numerator and denominator vanish, and is infinite when only the denominator vanishes. Let $\ud\muT=T^{-1}\ud t\ud\mu$ be the time-augmented probability measure on $(0,T)\times\RR^{2d}$, so that $$\Ent_{\muT}(g)= \frac{1}{T} \int_0^T\Ent_\mu(g_t)\ud t,$$ and define the following time-averaged distances, velocity Fisher information, and quadratic action:
\begin{equation} \label{eq:intro-averages}
      \begin{aligned}
 W_{x,T}^2(g)&=\frac1T\int_0^T W_2^2(q_t\mux,\mux)\ud t,
 && I_{v,T}(g)=\frac1T\int_0^T \Iv(g_t)\ud t, \\
 W_{v,T}^2(g)&=\frac1T\int_0^T\!\int q_t(x)
W_2^2(h_{t,x}\kv,\kv)\ud\mux(x)\ud t, &&J_T(g,W) = \int\frac{|W|^2}{g}\ud\muT.
\end{aligned}
\end{equation}
Here $q_t=\Piv g_t$, and $h_{t,x}=g_t(x,\cdot)/q_t(x)$ is the conditional velocity density at $x$, with $h_{t,x}=1$ on $\{q_t=0\}$. The functional $J_T(g, W)$ is the quadratic action of the density-current pair $(g, W)$ in \eqref{eq:action-intro}, and minimizing over the currents $W$ driving $g$ then gives the intrinsic action $\mathcal A_T(g)$. 

\begin{theorem}[Space-time HWI inequality]\label{thm:main-hwi}
Assume \eqref{eq:intro-potential}. For a given $T > 0$, let $(g,W)$ satisfy
\eqref{eq:controlled-intro} and \eqref{eq:intro-path}, with $g_t$ a probability density relative to $\mu$ at every time, $t\mapsto g_t\mu$ narrowly continuous on $[0,T]$, and $g_0\mu\in\PPtwo(\RR^{2d})$. Then $t\mapsto g_t\mu$ belongs to $C([0,T];\PPtwo(\RR^{2d}))$ with $\Ent_\mu(g_t)<\infty$ for every $t\in[0,T]$, and the Wasserstein current $\mathcal C_{\mux}(g_t)$ is continuous on $[0,T]$. Moreover, we have 
\begin{align}\label{eq:main-hwi}
 \Ent_{\muT}(g)\le{}&
 \frac{\mathcal C_{\mux}(g_0)-\mathcal C_{\mux}(g_T)}{T}
 +W_{x,T}(g)\sqrt{\mathcal A_T(g)}
 +3W_{v,T}(g)\sqrt{I_{v,T}(g)}\nonumber\\
 &-\frac K2W_{x,T}^2(g)-\frac32W_{v,T}^2(g).
\end{align}
\end{theorem}

For a static, velocity-independent density $g_t(x,v)=q(x)$ with $q\mux\in\PPtwo(\RR^d)$, finite entropy, and finite spatial Fisher information, we have $\mathcal A_T(g)=\int|\nabla_xq|^2/q\ud\mux$, and the space-time HWI \eqref{eq:main-hwi} reduces exactly to the standard HWI \eqref{eq:classical-hwi}; see \cref{rem:static} for more details. For time-dependent paths on a finite interval, the boundary term $(\mathcal C_{\mux}(g_0)-\mathcal C_{\mux}(g_T))/T$ cannot simply be omitted from \eqref{eq:main-hwi}; see \cref{rem:endpoint}. Under a compatible  Talagrand's transport inequality, \cref{thm:long-lsi} gives a boundary-free space-time LSI for every $T>0$, with a coefficient of the velocity Fisher information that grows as $T^{-2}$ at short times. This time dependence is optimal by \cref{rem:lsi-short-time}.

More generally, the spatial comparison measure $\mux$ in \cref{thm:main-hwi} can be replaced by any $\sigma\in\PPtwo(\RR^d)$ with $\Ent(\sigma\mid\mux)<\infty$. Let $T_{q,\sigma}$ be the optimal transportation map from $q\mux$ to $\sigma$, and define
\begin{equation}\label{eq:observables-intro}
 \Ds(g) := \frac12W_2^2(q\mux,\sigma),
 \qquad
  \Cs(g) := \int j(x)\cdot\bigl(x-T_{q,\sigma}(x)\bigr)\ud\mux.
\end{equation}
Here $\Cs$ is the Wasserstein current relative to the comparison measure $\sigma$, and $\xi_{q,\sigma}(x)=x-T_{q,\sigma}(x)$ is the Brenier displacement. For paths satisfying the assumptions of \cref{thm:main-hwi}, the Wasserstein first-variation formula gives $\frac{\ud}{\ud t}\Ds(g_t)=\Cs(g_t)$.
Clearly, the case $\sigma=\mux$ gives $\mathcal D_{\mux}(g)=\frac12W_2^2(q\mux,\mux)$ and the current $\mathcal C_{\mux}$ in \eqref{eq:cot-intro}. Similarly to \eqref{eq:intro-averages}, we write $$W_{x,\sigma,T}^2(g)=\frac1T\int_0^TW_2^2(q_t\mux,\sigma)\ud t.$$ Under the assumptions of \cref{thm:main-hwi}, the current $\Cs(g_t)$ is continuous on $[0,T]$, and we have the general space-time HWI:
\begin{equation} \label{eq:intro-comparator-hwi}
      \begin{aligned}
 \Ent_{\muT}(g)-\Ent(\sigma\mid\mux)\le{}&
 \frac{\Cs(g_0)-\Cs(g_T)}{T}
 +W_{x,\sigma,T}(g)\sqrt{\mathcal A_T(g)}
 +3W_{v,T}(g)\sqrt{I_{v,T}(g)} \\
 &-\frac K2W_{x,\sigma,T}^2(g)-\frac32W_{v,T}^2(g),
\end{aligned}
\end{equation}
see \cref{thm:comparator-hwi}. Noting $\mu=\mux\otimes\kv$, we have $\Ent(\sigma\otimes\kv\mid\mu)=\Ent(\sigma\mid\mux)$, and then the left side of \eqref{eq:intro-comparator-hwi} compares entropies with respect to the same measure $\mu$,
with $\sigma\otimes\kv$ as the static comparator. 
 
Our second main result is the kinetic analogue of the evolution variational inequality \eqref{eq:overdamped-evi-intro}.
We define the conditional velocity entropy by
\begin{equation} \label{eq:velocityentropy}
      \Entv(g)=\int q(x)\Ent_\kv(h_x)\ud\mux(x),
\end{equation}
where $q=\Piv g$ and $h_x=g(x,\cdot)/q(x)$ is the conditional velocity density at $x$, with $h_x=1$ on $\{q=0\}$. The chain rule for relative entropy under disintegration in $x$ gives $\Ent_\mu(g)=\Ent_{\mux}(q)+\Entv(g)$.

\begin{theorem}[Kinetic evolution variational inequality]\label{thm:kinetic-evi}
Assume \eqref{eq:intro-potential}. Let $g$ be a probability density relative to $\mu$ with $\Ent_\mu(g)<\infty$ and $g\mu\in\PPtwo(\RR^{2d})$, and set $g_t=\Pt_tg$ as its Langevin path \eqref{eq:langevinpath}. Fix $\sigma\in\PPtwo(\RR^d)$ with $\Ent(\sigma\mid\mux)<\infty$, and let $\Ds$ and $\Cs$ be as in \eqref{eq:observables-intro}. Then the function $t\mapsto\Ds(g_t)$ is continuously differentiable on $[0,\infty)$ with derivative $\Cs(g_t)$. Moreover, in $\mathcal D'(0,\infty)$,
\begin{equation}\label{eq:intro-langevin-comparison}
 \begin{aligned}
 &\frac{\ud^2}{\ud t^2}\Ds(g_t)
 +\gamma\frac{\ud}{\ud t}\Ds(g_t)
 +\Ent_\mu(g_t)-\Ent(\sigma\mid\mux)+K\Ds(g_t)\\
 &\qquad\le3\Entv(g_t)\le\frac32 \Iv(g_t).
 \end{aligned}
\end{equation}
\end{theorem}

By \cref{prop:semigroup-paths}, the entropy dissipation identity \eqref{eq:entropy-dissipation-intro} holds in time-integrated form along $g_t=\Pt_tg$, and hence the entropy and velocity Fisher information terms in \eqref{eq:intro-langevin-comparison} are locally integrable in time.

\subsection*{Proof ideas}
Both Theorems~\ref{thm:main-hwi} and \ref{thm:kinetic-evi} follow from a differential inequality for $\Cs(g_t)$ along controlled kinetic paths: integration in time yields the space-time HWI inequality, while specialization to Langevin evolution yields the kinetic EVI. The boundary term in the HWI retains information about the spatial current at both endpoints, as expected for an inertial evolution.

To explain the idea, we first compare gradient flow with inertial dynamics in finite-dimensional optimization. Fix a comparison point $y\in\RR^d$. Along the gradient flow $\dot x=-\nabla U(x)$, the Hessian lower bound \eqref{eq:intro-potential} directly gives the finite-dimensional EVI (a dynamical expression of the $K$-convexity of $U$):
\begin{equation}\label{eq:gradient-flow-intro}
 \frac{\ud}{\ud t}\frac12|x-y|^2 +U(x)-U(y)+\frac K2|x-y|^2\le0.
\end{equation}
This is the exact analogue of \eqref{eq:overdamped-evi-intro}, with Wasserstein distance replaced by Euclidean distance.

The finite-dimensional analogue of the underdamped Langevin dynamics is the heavy-ball equation~\cite{Polyak1964}, obtained by removing the noise from the Langevin equation,
\begin{equation}\label{eq:heavy-ball-ode-intro}
 \ddot x+\gamma\dot x+\nabla U(x)=0.
\end{equation}
Along its solutions, the first time derivative of the squared distance is
\begin{equation}\label{eq:heavy-ball-first-derivative-intro}
 \frac{\ud}{\ud t}\frac12|x-y|^2=\dot x\cdot(x-y),
\end{equation}
which does not directly give information on $U$.
Only the second time derivative yields a bound on the difference of objectives $U(x)-U(y)$, by using the heavy-ball equation \eqref{eq:heavy-ball-ode-intro} and the $K$-convexity of $U$,
\begin{equation}\label{eq:heavy-ball-intro}
 \frac{\ud}{\ud t}\bigl[\dot x\cdot(x-y)\bigr]
 +U(x)-U(y)+\frac K2|x-y|^2
 \le|\dot x|^2-\gamma\,\dot x\cdot(x-y),
\end{equation}
which is the inertial analogue of \eqref{eq:gradient-flow-intro}. 

For a controlled kinetic path, the spatial marginal $q_t\mux$ takes the place of $x(t)$, and the fixed spatial measure $\sigma$ takes the place of $y$. The spatial continuity equation gives the evolution of the squared Wasserstein distance (analogue of $|x-y|^2$): $\frac{\ud}{\ud t}\Ds(g_t)=\Cs(g_t)$, so the Wasserstein current plays the role of $\dot x\cdot(x-y)$ in \eqref{eq:heavy-ball-first-derivative-intro}.
This suggests differentiating $\Cs(g_t)$ in time and using the
$K$-geodesic convexity of the spatial relative entropy to obtain
an \emph{entropy comparison inequality}:
\begin{equation}\label{eq:intro-comparison}
 \frac{\ud}{\ud t}\Cs(g_t)
 +\Ent_\mu(g_t)-\Ent(\sigma\mid\mux)+K\Ds(g_t)
 \le3\Entv(g_t)+\int W_t\cdot\bigl(x-T_{q_t,\sigma}(x)\bigr)\ud\mu,
\end{equation}
which is analogous to \eqref{eq:heavy-ball-intro}. Here $W_t$ is the velocity current \eqref{eq:controlled-intro}. 
The relative entropy difference in \eqref{eq:intro-comparison} plays the role of $U(x)-U(y)$ in \eqref{eq:heavy-ball-intro}, and the Hessian lower bound yields the curvature term $K\Ds(g_t)$. The term $3\Entv(g_t)$ takes the place of the kinetic energy $|\dot x|^2$. The conditional velocity entropy vanishes when the conditional velocity laws equal the Gaussian equilibrium; by the Gaussian LSI, it can be bounded by the velocity Fisher information. For Langevin evolution, the velocity current is $W_t=-\gamma\nabla_vg_t$, and the last control term in \eqref{eq:intro-comparison} equals $-\gamma\Cs(g_t)$, the analogue of $-\gamma\,\dot x\cdot(x-y)$.

\medskip

For smooth paths with regular optimal transportation maps, differentiating $\Cs(g_t)$ with respect to time formally gives an exact identity whose nonnegative remainders can be discarded to obtain \eqref{eq:intro-comparison}; see Section~\ref{sec:identity}, which also describes the equality cases. The proof of the finite-action comparison inequality, \cref{thm:weak}, is independent of this formal identity.

In particular, for the rigorous proof, we will combine the Wasserstein acceleration inequality of~\cite{LiLu2026}*{Lemma 2.7} with integration by parts for the distributional Jacobian of the optimal transportation map. We first prove the entropy chain rule by localization, and then construct smooth approximations with convergence of the quadratic action and the total and spatial marginal entropies; see Appendix~\ref{app:recovery}. This yields \eqref{eq:intro-comparison} in distributions.

\subsection*{Consequences and comparison with our earlier work}
As mentioned above, the space-time HWI inequality gives a space-time LSI under positive curvature or under semiconvexity and a compatible  Talagrand's transport inequality; see \cref{cor:boundary-lsi}. 

In the rest of this subsection, we assume \eqref{eq:intro-potential} with $K\le0$, the  Talagrand's transport inequality with constant $\rho$, and $\alpha:=1+K/\rho>0$. Then \cref{thm:long-lsi} gives a universal constant $C$ such that, for every path satisfying the assumptions of \cref{thm:main-hwi} and every $T>0$,
\begin{equation}\label{eq:intro-boundary-free-lsi}
 \Ent_{\muT}(g)\le\frac{C}{\alpha^2}
       \left[\max\left\{1,\frac1{\rho T^2}\right\}I_{v,T}(g)
                         +\rho^{-1}\mathcal A_T(g)\right].
\end{equation}
The space-time LSI in our previous work~\cite{LiLu2026} assumes convexity and a spatial LSI. The estimate \eqref{eq:intro-boundary-free-lsi} extends it to semiconvex potentials under a compatible  Talagrand's transport inequality. On a more technical level, it also removes the additional growth and moment assumptions imposed in~\cite{LiLu2026}; see Section~\ref{sec:inequalities} for the precise comparison.

For Langevin evolution, the entropy estimate in \cref{prop:entropy-drop} quantifies relaxation from the initial time onward. \Cref{cor:decay} gives a universal constant $c>0$ such that, for every $\gamma>0$, every $t\ge0$, and every probability density $g$ of finite relative entropy,
\begin{equation}\label{eq:intro-decay}
 \Ent_\mu(\Pt_tg)
 \le\exp\!\left[-c\,\alpha^2
       \min\{\rho\gamma t^3,\gamma t,\rho t/\gamma\}\right]\Ent_\mu(g).
\end{equation}
The short-time powers in \eqref{eq:intro-boundary-free-lsi} and \eqref{eq:intro-decay} reflect the delayed dissipation of spatial nonequilibrium: a density depending only on $x$ is already at velocity equilibrium, so Hamiltonian transport must first make its spatial variation visible to velocity dissipation. Spatial shifts of Gaussian equilibria show that the inverse-square factor $O(T^{-2})$ in the space-time LSI and the cubic entropy-decay power $O(t^3)$ are both optimal; see \cref{rem:lsi-short-time,rem:cubic-sharp}.

At long times, the rate is proportional to $\alpha^2\min\{\gamma,\rho/\gamma\}$. For convex potentials, this recovers the rate order of~\cite{Lu2026}, which assumes a spatial LSI. The universal constants in \eqref{eq:intro-boundary-free-lsi} and \eqref{eq:intro-decay} do not depend on the dimension, the friction coefficient, or the potential, which enters only through $\rho$ and $\alpha$.

\subsection*{Geometric interpretation and related work}
  In Villani's approach to hypocoercivity, correction terms coupling spatial and velocity gradients recover spatial gradient control from velocity dissipation~\cite{Villani2009}. Dolbeault, Mouhot, and Schmeiser~\cite{DolbeaultMouhotSchmeiser2015} obtain decay through a modified $L^2$ norm that couples the spatial density to the nonequilibrium velocity component. The Wasserstein current serves a similar purpose for entropy: differentiating $\Cs(g_t)$ with respect to time gives spatial entropy control through the Brenier displacement, as in \eqref{eq:intro-comparison}, without requiring spatial Fisher information.

Geometrically, a particle at position $x$ with velocity $v$ contributes $v\cdot\xi_{q,\sigma}(x)$ to $\Cs(g)$, which is negative when the particle moves toward its target $T_{q,\sigma}(x)$. Changing only the velocities leaves the spatial marginal and the optimal transportation map fixed, but changes the Wasserstein current. Hence, velocity dissipation, although it acts only on velocities, can directly influence $\Cs(g_t)$, and \eqref{eq:intro-comparison} translates this change into control of the spatial entropy.
 
Space-time functional inequalities provide an alternative route for recovering coercivity by integrating over a time interval. Albritton, Armstrong, Mourrat, and Novack~\cite{AlbrittonArmstrongMourratNovack2024} developed kinetic Poincar\'e inequalities that combine velocity derivatives with the transport operator $\partial_t+\La$. Cao, Lu, and Wang~\cite{CaoLuWang2023} used this approach to obtain explicit $L^2$ convergence rates for underdamped Langevin dynamics. The space-time LSI of~\cite{LiLu2026} and the HWI proved here give corresponding entropy estimates.

The connection between transport and entropy also extends to the Schr\"odinger problem. Gentil, L\'eonard, Ripani, and Tamanini~\cite{GentilLeonardRipaniTamanini2018} prove the classical HWI through entropic interpolations and a zero-noise limit. For Langevin reference dynamics, Chiarini, Conforti, Greco, and Ren~\cite{ChiariniConfortiGrecoRen2022} establish long-time entropy estimates for kinetic Schr\"odinger bridges, which minimize relative entropy on path space subject to prescribed endpoint marginals. Our space-time HWI applies to every controlled path satisfying the assumptions of \cref{thm:main-hwi}, without requiring the path to minimize a cost.

Hamiltonian transport has a Lie--Poisson formulation for Vlasov equations~\cites{Morrison1980,MarsdenWeinstein1982,MillerEtAl2022}, and its combination with symmetric velocity dissipation appears in variational descriptions of Vlasov--Fokker--Planck equations~\cite{DuongPeletierZimmer2013}. In  \cref{rem:structure}, we express the Wasserstein current (i.e., the damping term) in \eqref{eq:intro-langevin-comparison} through this structure. Our intrinsic action measures additional acceleration relative to Hamiltonian dynamics. Brigati, Maas, and Quattrocchi~\cite{BrigatiMaasQuattrocchi2025} minimize squared acceleration over paths with prescribed endpoints. In $\mathcal A_T(g)$, the path is fixed, and only the current driving it is optimized; see \cref{rem:otikin} for the detailed discussion.

\subsection*{Organization} In Section~\ref{sec:setting}, we fix the notation and recall the elementary entropy and action estimates. Section~\ref{sec:identity} formally derives the exact identity for the evolution of the Wasserstein current and deduces the semiconvex comparison inequality. These computations illustrate the main ideas behind the rigorous argument. 
In Section~\ref{sec:weak}, we establish the comparison inequality for finite-action paths under semiconvexity and obtain both main theorems: \cref{thm:main-hwi} by integrating in time, and \cref{thm:kinetic-evi} by specializing to Langevin paths. As a consequence, Section~\ref{sec:inequalities} proves the LSI with a nonpositive boundary term and the boundary-free LSI for every $T>0$, and then Section~\ref{sec:decay} proves entropy decay on each time interval under semiconvexity, and the resulting cubic short-time bound and exponential decay for every $\gamma>0$. Appendix~\ref{app:recovery} proves the entropy chain rule, verifies its application to Langevin paths, and constructs the approximation used in \cref{thm:weak}.

\subsection*{AI Usage} The main mathematical ideas in this work are from the authors. We used LLM-based AI assistants in editing and polishing the manuscript. The authors take full responsibility for the manuscript.

\section{Setting and preliminaries}\label{sec:setting}

Throughout, we work on the phase space $\RR^{2d}$ with the reference Gibbs measure $\mu=\mux\otimes\kv$ from \eqref{eq:gibbs-intro}. We assume that $U\in C^\infty(\RR^d)$, $Z_x=\int\e^{-U}\ud x<\infty$, and $\mu\in\PPtwo(\RR^{2d})$. For a probability density $g$ relative to $\mu$, we write $\Ent_\mu(g)=\int g\log g\ud\mu$. Equivalently, for $\nu\ll\mu$, we write $\Ent(\nu\mid\mu)=\Ent_\mu(\ud\nu/\ud\mu)$, and use the corresponding notation for entropy relative to the spatial and velocity marginals $\mux$ and $\kv$. When convenient, we identify an absolutely continuous measure with its Lebesgue density; for example, $q\mux$ also denotes $q\e^{-U}/Z_x$, and $\sigma/\mux$ is the density of $\sigma$ relative to $\mux$. For an open set $\Omega$, $\mathcal D'(\Omega)$ denotes the space of distributions on $\Omega$.

The operators $\La$ and $\Ls$ are defined in \eqref{eq:generators-intro}. On smooth compactly supported functions, $\La$ is skew-symmetric in $L^2(\mu)$, while $\Ls=-\nabla_v^*\nabla_v$ is symmetric and nonpositive. We use the formal adjoints
\[
 \nabla_v^*F=-\divg_vF+v\cdot F,
 \qquad \nabla_x^*F=-\divg_xF+\nabla U\cdot F,
\]
with $\nabla_x^*$ acting row by row on matrix fields. 

\subsection{Conditional velocity laws and entropy}

For a probability density $g$ relative to $\mu$ with $g\mu\in\PPtwo(\RR^{2d})$, let $\Piv$ denote integration with respect to $\ud \kv(v)$.
Following~\cites{Lu2026,LiLu2026}, we define the spatial marginal density $q$, the spatial current $j$, the conditional velocity density $h_x$, and the conditional mean velocity $m$ by
\begin{equation}\label{eq:moments}
 q=\Piv g,\qquad j=\Piv(vg),\qquad  h_x=g(x,\cdot)/q(x),\qquad m=j/q.
\end{equation}
The ratios defining $h_x$ and $m$ are understood on $\{q>0\}$ with the convention $h_x=1$ and $m=j=0$ on $\{q=0\}$. The conditional covariance is 
\begin{equation}\label{eq:covariance}
 \Sigma_g(x)=\int(v-m(x))\otimes(v-m(x))h_x(v)\ud\kv(v).
\end{equation}
Disintegration with respect to $x$ also factorizes the entropy: 
\begin{equation}\label{eq:conditional-entropy}
 \Entv(g)=\int q(x)\Ent_\kv(h_x)\ud\mux(x),
 \qquad
 \Ent_\mu(g)=\Ent_{\mux}(q)+\Entv(g).
\end{equation}
Here $\Entv(g)$ is the conditional velocity entropy, and \eqref{eq:conditional-entropy} holds in the extended sense, with values in $[0,\infty]$. We also recall the velocity Fisher information, for a phase-space density $g$, 
\begin{equation}\label{eq:fisher-definition}
 \Iv(g)=4\int|\nabla_v\sqrt g|^2\ud\mu,
\end{equation}
which agrees with the expression $\int|\nabla_vg|^2/g\ud\mu$ in \eqref{eq:entropy-dissipation-intro} when $g$ is smooth and strictly positive. For general densities, the gradient in \eqref{eq:fisher-definition} is understood weakly, and $\Iv(g) :=+\infty$ if it does not exist. When $\Iv(g)<\infty$, the velocity score $\nabla_v\log g$ is defined as $2\nabla_v\sqrt g/\sqrt g$ on $\{g>0\}$ and as zero elsewhere.

For $g\mu\in\PPtwo(\RR^{2d})$ with $\Ent_\mu(g)<\infty$, comparison with the Gaussian $\mathcal N(m(x),\Id)$ gives
\begin{equation*}
       \Ent_\kv(h_x) = \frac12|m(x)|^2
  +\Ent\bigl(h_x\kv\mid\mathcal N(m(x),\Id)\bigr) \ge\frac12|m(x)|^2.
\end{equation*}
Integrating against $q\mux$ and applying the Gaussian LSI to the conditional velocity laws yields
\begin{equation}\label{eq:mean-entropy}
      \int q|m|^2\ud\mux\le2\Entv(g) \le \Iv(g).
\end{equation}

The conditional mean velocity $m(x)$ is also related directly to the velocity score $\nabla_v \log g$. Indeed, if $\Iv(g)<\infty$, Gaussian integration by parts gives
\begin{equation}\label{eq:score-mean}
 \int\nabla_v\log g(x,v)h_x(v)\ud\kv(v)=m(x),
 \qquad q\mux\text{-a.e. }x.
\end{equation}
Setting $W_v^2(g)=\int q(x)W_2^2(h_x\kv,\kv)\ud\mux(x)$, the Gaussian HWI inequality~\cite{OttoVillani2000} gives \begin{equation}\label{eq:fiber-hwi}
 \Entv(g)\le W_v(g)\sqrt{\Iv(g)}-\tfrac12W_v^2(g).
\end{equation}

To identify the contributions of the mean and covariance to velocity entropy, we compare a probability law with the Gaussian having the same mean and covariance.
\begin{lemma}[Gaussian entropy decomposition]\label{lem:gaussian}
Let $\eta\in\PPtwo(\RR^d)$ have finite entropy relative to $\kv$, mean $m$, and covariance $\Sigma$. Then $\Sigma$ is positive definite, and
\begin{equation}\label{eq:gaussian-exact}
 2\Ent(\eta\mid\kv) - 2\Ent\bigl(\eta\mid\mathcal N(m,\Sigma)\bigr)
 =|m|^2+\tr\Sigma-d-\log\det\Sigma.
\end{equation}
\end{lemma}

\begin{proof}
Finite entropy relative to $\kv$ implies absolute continuity of $\eta$, and hence its covariance $\Sigma$ is positive definite. The Gaussian density ratio satisfies
\[
 \log\frac{\ud\mathcal N(m,\Sigma)}{\ud\kv}(v)
 =\frac12\bigl[|v|^2-(v-m)\cdot\Sigma^{-1}(v-m)-\log\det\Sigma\bigr],
\]
which is $\eta$-integrable. Since $\eta$ has mean $m$ and covariance $\Sigma$,
\[
 \begin{aligned}
 \Ent(\eta\mid\kv)-\Ent\bigl(\eta\mid\mathcal N(m,\Sigma)\bigr)
 &=\int\log\frac{\ud\mathcal N(m,\Sigma)}{\ud\kv}\ud\eta\\
 &=\frac12\bigl[|m|^2+\tr\Sigma-\tr(\Sigma^{-1}\Sigma)-\log\det\Sigma\bigr],
 \end{aligned}
\]
which proves \eqref{eq:gaussian-exact}.
\end{proof}

\subsection{Bounds on the Wasserstein current}

Let $g\mu\in\PPtwo(\RR^{2d})$ and fix $\sigma\in\PPtwo(\RR^d)$ with $\Ent(\sigma\mid\mux)<\infty$. Brenier's theorem~\cites{brenier1991polar,villani2021topics} gives a $q\mux$-a.e. unique optimal transportation map $T_{q,\sigma}$ from $q\mux$ to $\sigma$. We recall the Brenier displacement
\[\xi_{q,\sigma}=x-T_{q,\sigma}(x).\]
We recall the definitions of $\Ds$ and $\Cs$ from \eqref{eq:observables-intro}. For $\sigma=\mux$, we write $T_q = T_{q,\mux}$, $\xi_q = \xi_{q,\mux}$, and $\mathcal D_{\mux}(g) = \frac12 W_2^2(q\mux,\mux)$, with corresponding current $\mathcal C_{\mux}$. 
These functionals also extend to phase-space laws
$\nu\in\PPtwo(\RR^{2d})$ with spatial marginal $q\mux$:
\[
 \Ds(\nu)=\frac12W_2^2(q\mux,\sigma),
 \qquad
 \Cs(\nu)=\int v\cdot\xi_{q,\sigma}(x)\ud\nu(x,v),
\]
which agree with \eqref{eq:observables-intro} for $\nu=g\mu$.

By Cauchy--Schwarz and \eqref{eq:score-mean}, we have 
\begin{equation}\label{eq:current-size}
 |\Cs(g)|^2\le W_2^2(q\mux,\sigma)\int q|m|^2\ud\mux
                 \le W_2^2(q\mux,\sigma)\Iv(g).
\end{equation}
To control the current by entropy, suppose that $\mux$ satisfies the Talagrand's transportation inequality \cites{Talagrand1996,OttoVillani2000,BakryGentilLedoux2014}: for every probability density $q$ relative to $\mux$, 
\begin{equation}\label{eq:spatial-t2}
 W_2^2(q\mux,\mux)\le\frac2\rho\Ent_{\mux}(q),\qquad \rho>0.
\end{equation}
 Combining \eqref{eq:current-size} with \eqref{eq:mean-entropy} and \eqref{eq:spatial-t2} gives
\begin{equation}\label{eq:current-t2}
 |\mathcal C_{\mux}(g)|^2 \le\frac2\rho\Ent_{\mux}(q)\Iv(g)
, \qquad
 |\mathcal C_{\mux}(g)|\le\rho^{-1/2}\Ent_\mu(g),
\end{equation}
which are the corrector estimates of~\cite{Lu2026}*{Lemma 3.2}. 
Under \eqref{eq:spatial-t2}, finite entropy also implies a finite second moment. Indeed, the triangle inequality for $W_2$, followed by the spatial and conditional Gaussian Talagrand inequalities, gives
\begin{align}\label{eq:entropy-second-moment}
 \int q|x|^2\ud\mux
   &\le\frac4\rho\Ent_{\mux}(q)+2\int|x|^2\ud\mux,\nonumber\\
 \int g|v|^2\ud\mu&\le4\Entv(g)+2d.
\end{align}

\subsection{Controlled paths and the action}

Following~\cite{LiLu2026}*{Section 3.2}, let $g=\{g_t\}_{t\in[0,T]}$ be a path of probability densities with respect to $\mu$. A measurable velocity current $W:(0,T)\times\RR^{2d}\to\RR^d$ drives $g$ if \eqref{eq:controlled-intro} holds in $\mathcal D'((0,T)\times\RR^{2d})$. The quadratic action of the density-current pair $(g,W)$ is $J_T(g,W)$ from \eqref{eq:intro-averages}.

As in~\cite{LiLu2026}, let $\varrho_\infty$ denote the Lebesgue density of $\mu$, and write $\varrho_t=g_t\varrho_\infty$ and $F_t=W_t\varrho_\infty$ for the corresponding Lebesgue density and velocity current. We also introduce $\nu_t=g_t\mu=\varrho_t\ud x\ud v$ for the phase-space law.
Then \eqref{eq:controlled-intro} becomes
\begin{equation}\label{eq:ordinary-control}
 \partial_t\varrho+\divg_x(v\varrho)+
 \divg_v\bigl(-\nabla U\,\varrho+F\bigr)=0.
\end{equation}
If $J_T(g,W)<\infty$, one may write $W=gA$ with a velocity field $A\in L^2(g\ud\muT;\RR^d)$. In this representation, $F=\varrho A$ and $J_T(g,W)=\int g|A|^2\ud\muT$.

The evolution of the zeroth and first velocity moments of $\varrho_t$ under 
\eqref{eq:controlled-intro}, as in~\cite{LiLu2026}*{proof of Lemma 3.5 and Lemma B.17}, satisfies
\begin{equation}\label{eq:moment-equations}
 \partial_tq=\nabla_x^*j,\qquad
 \partial_tj-\nabla_x^*\!\left(\frac{j\otimes j}{q}\right)
      =\nabla_x^*(q\Sigma_g)-q\nabla U+\Piv W,
\end{equation}
where $\Sigma_g$ is the conditional covariance in \eqref{eq:covariance}. The first equation is the continuity equation for the spatial law $q\mux$ with velocity $m=j/q$. The second involves the covariance stress $q\Sigma_g$, so the evolution of $q$ and $j$ is not closed; we retain the full phase-space law $g\mu$ as the state variable. If $J_T(g,W)<\infty$ and the second moments are locally bounded in time, velocity cutoffs justify \eqref{eq:moment-equations} in $\mathcal D'((0,T)\times\RR^d)$.

Minimizing the quadratic action over all currents driving $g$, that is, satisfying \eqref{eq:controlled-intro} in distributions, gives the intrinsic action of the path:
\begin{equation}\label{eq:action}
 \mathcal A_T(g)=\inf_{W\text{ drives }g}J_T(g,W).
\end{equation}
The infimum is $+\infty$ if no current drives $g$.

For smooth density-current pairs $(g,W)$ satisfying \eqref{eq:controlled-intro}, with $g=\delta+f$, where $\delta>0$ is a constant, and $f(x,v) \ge 0$ and $W(x,v)$ are compactly supported uniformly in time,  integration by parts gives the entropy chain rule, as in~\cite{LiLu2026}*{Lemma 3.4},
\begin{equation}\label{eq:entropy-chain}
 \Ent_\mu(g_t)-\Ent_\mu(g_s)
    =\int_s^t\!\int W_r\cdot\nabla_v\log g_r\ud\mu\ud r.
\end{equation} 
Under the finite-action assumptions of \cref{thm:main-hwi}, \cref{prop:entropy-chain} proves \eqref{eq:entropy-chain} by a cutoff argument for every $0\le s\le t\le T$, including both endpoints.

\begin{remark}[Relation to kinetic optimal transport]\label{rem:otikin}
In the representation $W=gA$, the field $A$ gives the additional acceleration relative to Hamiltonian transport. To compare the intrinsic action with kinetic optimal transport \cite{BrigatiMaasQuattrocchi2025}, fix a path $\nu_t=g_t\mu$ and let $b_t=-\nabla U+A_t$ be the total acceleration. Equation \eqref{eq:ordinary-control} and the action \eqref{eq:action} become
\[
 \partial_t\nu_t+v\cdot\nabla_x\nu_t+\divg_v(b_t\nu_t)=0,
 \qquad
 \mathcal A_T(g)=\inf_b\frac1T\int_0^T\!\int
                  |b_t+\nabla U|^2\ud\nu_t\ud t,
\]
where the infimum is over fields $b$ driving the prescribed path $\nu$. The kinetic Benamou--Brenier formula of~\cite{BrigatiMaasQuattrocchi2025}*{Theorem 4.10} uses the same continuity equation to characterize the fixed-horizon discrepancy $\tilde{\mathsf d}_T$:
\[
 \tilde{\mathsf d}_T^2(\nu_0,\nu_T)
       =\inf_{(\eta,b)}T\int_0^T\!\int|b_t|^2\ud\eta_t\ud t.
\]
Here the infimum is over paths $\eta\in C([0,T];\PPtwo(\RR^{2d}))$ and square-integrable acceleration fields $b$ satisfying the kinetic continuity equation, with $\eta_0=\nu_0$ and $\eta_T=\nu_T$. Thus $\tilde{\mathsf d}_T$ minimizes the total acceleration cost over paths joining two endpoints, whereas $\mathcal A_T(g)$ minimizes the cost of the additional acceleration $b+\nabla U$ along a prescribed path. In particular, Hamiltonian motion has zero intrinsic action but, in general, a nonzero total acceleration cost. We emphasize that our HWI inequality \eqref{eq:main-hwi} uses the usual $W_2$ distances for the spatial marginal and conditional velocity laws, rather than $\tilde{\mathsf d}_T$.
\end{remark}

\section{Evolution of the Wasserstein current}\label{sec:identity}

In this section, we explain the geometric meaning of the comparison inequality \eqref{eq:intro-comparison} and its connection to the evolution of the Wasserstein current. 
The computations in this section are formal, to keep the idea transparent.
We assume that the density path, the comparison measure $\sigma$, and the optimal transportation maps $T_t:=T_{q_t,\sigma}$ are smooth enough, with enough decay, so that every differentiation and integration by parts below is valid. 
Moreover, suppose that each $T_t$ is a $C^2$ diffeomorphism of $\RR^d$ and the conditional covariance $\Sigma_{g_t}$ is invertible.

We recall that along a controlled path satisfying \eqref{eq:controlled-intro}, the spatial continuity equation \eqref{eq:moment-equations} and the Wasserstein first-variation formula~\cite{ambrosio2005gradient}*{Theorem 8.4.7} give the Wasserstein current as the time derivative of $\Ds(g_t) = \frac12W_2^2(q_t\mux,\sigma)$:
\begin{equation}\label{eq:controlled-observable}
 \frac{\ud}{\ud t}\Ds(g_t)=\Cs(g_t).
\end{equation}
We next derive the time derivative of $\Cs(g_t)$.

\begin{proposition}[Second-derivative identity]\label{thm:exact}
Along the controlled path,
\begin{equation}\label{eq:exact}
 \frac{\ud}{\ud t}\Cs(g_t)
 +\Ent_\mu(g_t)-\Ent(\sigma\mid\mux)
 +\mathcal R_U(q_t,\sigma)+\mathcal E_\sigma(t)
 =3\Entv(g_t)+\int W_t\cdot\xi_{q_t,\sigma}\ud\mu,
\end{equation}
where the potential remainder $\mathcal R_U(q_t,\sigma)$ and the nonnegative remainder $\mathcal E_\sigma(t)$ are given by \eqref{eq:potential-remainder} and \eqref{eq:exact-error} below, respectively. 
\end{proposition}

In \eqref{eq:exact}, the potential remainder is
\begin{equation}\label{eq:potential-remainder}
 \mathcal R_U(q,\sigma)=\int q(x)
 \bigl[U(T_{q,\sigma}(x))-U(x)
              +\nabla U(x)\cdot\xi_{q,\sigma}(x)\bigr]\ud\mux(x).
\end{equation}
A lower Hessian bound on $U$ is needed only to bound this term from below. With $G_t:=DT_t$ the Jacobian matrix of $T_t$, the remainder $\mathcal E_\sigma(t)$ is
\begin{equation}\label{eq:exact-error}
 \begin{aligned}
 \mathcal E_\sigma(t)={}&
 \int q_t\,
    (\partial_tT_t+G_tm_t)\cdot G_t^{-1}
                    (\partial_tT_t+G_tm_t)\ud\mux\\
 &+\int q_t\left[\tr(\Sigma_{g_t}G_t)-d
                           -\log\det(\Sigma_{g_t}G_t)\right]\ud\mux\\
 &+2\int q_t(x)\Ent\left(h_{t,x}\kv\,
                \middle|\,\mathcal N(m_t(x),\Sigma_{g_t}(x))\right)\ud\mux(x).
 \end{aligned}
\end{equation}
Each integral in \eqref{eq:exact-error} is nonnegative, since the Jacobian matrix $G_t$ of the optimal transportation map $T_t$ is symmetric positive definite. For the second integral, although $\Sigma_{g_t}G_t$ need not be symmetric, it is similar to the positive definite symmetric matrix $\Sigma_{g_t}^{1/2}G_t\Sigma_{g_t}^{1/2}$. Its eigenvalues are therefore positive, and $r-1-\log r\ge0$ for $r>0$ gives the claim. Thus, $\mathcal{E}_\sigma(t)$ can be discarded from \eqref{eq:exact}, yielding the desired comparison inequality \eqref{eq:smooth-comparison}. \cref{rem:nonnegative-terms} describes when each integral in \eqref{eq:exact-error} vanishes.

\begin{proof}[Proof of \cref{thm:exact}]
For simplicity, we often suppress the time subscript and write $T=T_t$, $G=G_t$, $\xi=x-T$, $m=m_t$, and $\Sigma=\Sigma_{g_t}$. Following the finite-dimensional heavy-ball computation leading to \eqref{eq:heavy-ball-intro}, we differentiate $\Cs(g_t)$ in time and identify the resulting entropy terms. 

\smallskip

\noindent \emph{Step 1: Differentiation along the Hamiltonian flow.} For a function $\psi_t(x,v)$, the controlled equation \eqref{eq:controlled-intro} and the skew-symmetry of $\La$ give
\begin{equation}\label{eq:weak-control}
 \frac{\ud}{\ud t}\int\psi_tg_t\ud\mu
 =\int g_t(\partial_t+\La)\psi_t\ud\mu
  +\int W_t\cdot\nabla_v\psi_t\ud\mu.
\end{equation}
The Wasserstein current is the integral of
$\psi_t(x,v)=v\cdot\xi(x)$ against $g_t\mu$. For this test function,
\begin{equation} \label{eq:derivative}
    \nabla_v\psi_t=\xi,\qquad
 \partial_t\psi_t=-v\cdot\partial_tT,\qquad
 \nabla_x\psi_t=(\Id-G)v,
\end{equation}
where the last identity uses the symmetry of the Jacobian
$\Id-G$ of $\xi$. 
Along a Hamiltonian flow $\dot x_t=v_t$, $\dot v_t=-\nabla U(x_t)$, the chain rule gives the change rate of $\psi_t$:
\[
 \frac{\ud}{\ud t}\psi_t(x_t,v_t)
 =\bigl[\partial_t\psi_t+v\cdot\nabla_x\psi_t
       -\nabla U\cdot\nabla_v\psi_t\bigr](x_t,v_t)
 =(\partial_t+\La)\psi_t(x_t,v_t). 
\]
Substituting the derivatives \eqref{eq:derivative} above yields
\begin{equation} \label{eq:expressham}
    (\partial_t+\La)\psi_t =-\nabla U\cdot\xi+|v|^2  -v\cdot(\partial_tT+Gv).
\end{equation}
As in the heavy-ball computation, the first two terms come from the force and the kinetic energy. The last term accounts for the motion of the comparison point $T_t(x_t)$ along the flow:
\[
 \frac{\ud}{\ud t}T_t(x_t) =\partial_tT_t(x_t)+G_t(x_t)v_t,
\]
unlike the fixed point $y$ in \eqref{eq:heavy-ball-intro}.

Using \eqref{eq:expressham}, we consider the integrals in \eqref{eq:weak-control} against $g\mu$, first in $v$ against the conditional velocity law $h_x\kv$, which has mean $m$ and second moment $\Sigma+m\otimes m$. Specifically, for every vector field $b(x)$ and matrix field $M(x)$, we have 
\[
 \int g\,v\cdot b\ud\mu=\int q\,m\cdot b\ud\mux,
 \qquad
 \int g\,v\cdot Mv\ud\mu=\int q\bigl[m\cdot Mm+\tr(\Sigma M)\bigr]\ud\mux.
\]
With $M=\Id$, the kinetic term becomes $\int g|v|^2\ud\mu=\int q(|m|^2+\tr\Sigma)\ud\mux$. With $b=\partial_tT$ and $M=G$, the term from the moving comparison point gives
\[
 \int g\,v\cdot(\partial_tT+Gv)\ud\mu
 =\int q\bigl[m\cdot(\partial_tT+Gm)+\tr(\Sigma G)\bigr]\ud\mux.
\]
The force term $\nabla U\cdot\xi$ does not depend on $v$, and the control term in \eqref{eq:weak-control} is $\int W\cdot\nabla_v\psi_t\ud\mu=\int W\cdot\xi\ud\mu$. Substituting these identities into \eqref{eq:weak-control} gives
\begin{equation}\label{eq:current-before-entropy}
 \begin{aligned}
 \frac{\ud}{\ud t}\Cs(g_t)
 ={}&\int q\bigl[|m|^2+\tr\Sigma-\tr(\Sigma G)\bigr]\ud\mux
       -\int q\nabla U\cdot\xi\ud\mux\\
 &-\int q\,m\cdot(\partial_tT+Gm)\ud\mux
       +\int W\cdot\xi\ud\mu.
 \end{aligned}
\end{equation}

\smallskip

\noindent
\emph{Step 2: Motion of the target points.}
In the term $-\int q\,m\cdot(\partial_tT+Gm)\ud\mux$, the vector field $\partial_tT+Gm$ describes the velocity of the target points along the flow of the conditional mean velocity $m$. Indeed, along a trajectory satisfying $\dot x_t=m_t(x_t)$,
\[
 \frac{\ud}{\ud t}T_t(x_t)  =\partial_tT_t(x_t)+G_t(x_t)m_t(x_t).
\]
The spatial continuity equation in \eqref{eq:moment-equations} shows that $m_t$ transports $q_t\mux$. Since $(T_t)_\#(q_t\mux)=\sigma$ and $\sigma$ is fixed, differentiating $\int f(T_t)q_t\ud\mux=\int f\ud\sigma$ gives, for every smooth test function $f$,
\[
 0=\frac{\ud}{\ud t}\int f(T_t)q_t\ud\mux =\int \nabla f(T)\cdot(\partial_tT+Gm)\,q\ud\mux.
\]
Thus the velocity $\partial_tT+Gm$ is orthogonal in $L^2(q\mux)$ to every
gradient field evaluated at $T$, i.e., $\nabla f\circ T$. 
Optimality of $T$ supplies such a gradient field. Writing $T_t=\nabla\phi_t$ for a convex Brenier potential $\phi_t$, we have $\partial_tT_t=\nabla\dot\phi_t$. For $f=\dot\phi_t\circ T_t^{-1}$, the chain rule and the symmetry of $G$ give $G\nabla f(T)=\partial_tT$, so that
\[
 \int q\,G^{-1}\partial_tT\cdot(\partial_tT+Gm)\ud\mux=0.
\]
Substituting $m=G^{-1}(\partial_tT+Gm)-G^{-1}\partial_tT$ yields
\begin{equation}\label{eq:map-orthogonality}
 \int q\,m\cdot(\partial_tT+Gm)\ud\mux
       =\int q(\partial_tT+Gm)\cdot G^{-1}
                        (\partial_tT+Gm)\ud\mux,
\end{equation}
which is the first integral in \eqref{eq:exact-error}.

\smallskip

\noindent
\emph{Step 3: Spatial entropy.} Since $T_\#(q\mux)=\sigma$, the change-of-variables formula reads
\[
 q\mux=\sigma(T)\det G.
\]
Dividing by $\mux$, taking logarithms, and integrating against $q\mux$, with $\int q\log\frac{\sigma}{\mux}(T)\ud\mux=\Ent(\sigma\mid\mux)$, gives
\begin{equation}\label{eq:transport-entropy}
 \Ent_{\mux}(q)-\Ent(\sigma\mid\mux)
   =\int q\bigl[\log\det G+U(x)-U(T(x))\bigr]\ud\mux.
\end{equation}
After adding $\mathcal R_U(q,\sigma)$ to both sides of \eqref{eq:transport-entropy}, the right side becomes $\int q(\log\det G+\nabla U\cdot\xi)\ud\mux$. Its force term cancels $-\int q\nabla U\cdot\xi\ud\mux$ in \eqref{eq:current-before-entropy}.

\smallskip

\noindent
\emph{Step 4: Covariance and Jacobian.} To combine the covariance and Jacobian terms, we use the matrix identity
\begin{equation}\label{eq:matrix-cancellation}
 \begin{aligned}
 \tr\Sigma-\tr(\Sigma G)+\log\det G
 ={}&\tr\Sigma-d-\log\det\Sigma -\bigl[\tr(\Sigma G)-d-\log\det(\Sigma G)\bigr].
 \end{aligned}
\end{equation}
When $\Sigma=\Id$, the bracket in \eqref{eq:matrix-cancellation} reduces to $\tr G-d-\log\det G$, which is nonnegative and discarded in the analogous computation for the overdamped EVI \eqref{eq:overdamped-evi-intro}. In the kinetic setting, the conditional covariance $\Sigma$ replaces $\Id$. Combining \eqref{eq:current-before-entropy}, \eqref{eq:map-orthogonality}, \eqref{eq:transport-entropy}, and \eqref{eq:matrix-cancellation} leaves the term
$|m|^2+\tr\Sigma-d-\log\det\Sigma$, which can be expressed as $2\Ent_\kv(h_x)-2\Ent(h_x\kv\mid\mathcal N(m,\Sigma))$, by \eqref{eq:gaussian-exact}. We therefore have
\[
 \frac{\ud}{\ud t}\Cs(g_t)+\Ent_{\mux}(q)-\Ent(\sigma\mid\mux)
       +\mathcal{R}_U+\mathcal{E}_\sigma
     =2\Entv(g)+\int W\cdot\xi\,\ud\mu,
\]
which proves \eqref{eq:exact} upon adding $\Entv(g)$ to both sides. 
\end{proof}

In addition, assuming $\nabla^2U\succeq K\Id$, Taylor's formula gives
\[
 U(T)-U(x)+\nabla U(x)\cdot(x-T)\ge\tfrac K2|x-T|^2,
\]
which further implies  $\mathcal R_U\ge K\Ds$. Then, discarding the nonnegative remainder $\mathcal E_\sigma$ in \eqref{eq:exact} yields
\begin{equation}\label{eq:smooth-comparison}
 \frac{\ud^2}{\ud t^2}\Ds(g_t)
 +\Ent_\mu(g_t)-\Ent(\sigma\mid\mux)+K\Ds(g_t)
 \le3\Entv(g_t)+\int W_t\cdot\xi_{q_t,\sigma}\ud\mu.
\end{equation}
Inequality \eqref{eq:smooth-comparison} is the analogue of \eqref{eq:heavy-ball-intro}, with $3\Entv(g_t)$ in place of $|\dot x|^2$. Here, the coefficient $3$ is sharp; see \cref{rem:sharp-three}.

\begin{remark}[Conditions for equality]\label{rem:nonnegative-terms}
Equality holds in \eqref{eq:smooth-comparison} precisely when $\mathcal R_U(q_t,\sigma)=K\Ds(g_t)$ and the three integrals in \eqref{eq:exact-error} vanish. The first integral vanishes if and only if  $\partial_tT_t+G_tm_t=0$, meaning that the optimal transportation map to the fixed target is transported by the flow of the conditional mean velocity $m_t$. The second integral vanishes when $\Sigma_{g_t}=G_t^{-1}$, while the third one vanishes when the conditional velocity law is Gaussian.
\end{remark}

For the Langevin evolution, the velocity current is $W_t=-\gamma\nabla_v g_t$. The score identity \eqref{eq:score-mean}, together with \eqref{eq:controlled-observable}, gives
\[
 \int W_t\cdot\xi_{q_t,\sigma}\,\ud\mu
 =-\gamma\Cs(g_t)=-\gamma\frac{\ud}{\ud t}\Ds(g_t),
\]
and then the comparison inequality \eqref{eq:smooth-comparison} yields the kinetic EVI \eqref{eq:intro-langevin-comparison}. 

\begin{remark}[Hamiltonian and dissipative structure]\label{rem:structure}
This remark explains how the first variations of $\Ds$ and $\Cs$ give rise to the control and damping terms in \eqref{eq:smooth-comparison} and \eqref{eq:intro-langevin-comparison}.  Let $\mathcal H(g)=\Ent_\mu(g)$, and write $\mathcal F_g$ for the first variation at $g$ of a functional $\mathcal F$, with respect to the pairing in $L^2(\mu)$. The canonical Poisson bracket and the symmetric form in the velocity variables are
\begin{equation}\label{eq:bracket}
 \begin{aligned}
 \{\mathcal F,\mathcal G\}(g)
 &=\int\bigl(\nabla_x\mathcal F_g\cdot\nabla_v\mathcal G_g
             -\nabla_v\mathcal F_g\cdot\nabla_x\mathcal G_g\bigr)g\ud\mu,\\
 \Gamma_v(\mathcal F,\mathcal G)(g)
 &=\int\nabla_v\mathcal F_g\cdot\nabla_v\mathcal G_g\,g\ud\mu.
 \end{aligned}
\end{equation}
Here $\{\mathcal F,\mathcal G\}$ is the usual Lie--Poisson bracket, and $\Gamma_v$ is induced by the velocity carr\'e du champ~\cite{Baudoin2017}. The first variation of entropy is $\mathcal H_g=\log g+\text{constant}$, implying $\Gamma_v(\mathcal H,\mathcal H)=\Iv(g)$. In $\{\mathcal F,\mathcal H\}(g)$, integration by parts against $\mu$ moves the derivatives from $g$ to $\mathcal F_g$, and the mixed second derivatives of $\mathcal F_g$ cancel by symmetry. Hence $\{\mathcal F,\mathcal H\}(g)=\int g\La\mathcal F_g\ud\mu$. Thus the underdamped Langevin equation combines the Hamiltonian flow of $\mathcal H$ for this bracket with $\gamma$ times the gradient flow of $\mathcal H$ for the form $\Gamma_v$. This conservative--dissipative interpretation is established in~\cite{DuongPeletierZimmer2013}. 

Along the controlled path, the chain rule for a functional $\mathcal F$ and \eqref{eq:weak-control} with $\psi=\mathcal F_{g_t}$ show that 
\[
 \frac{\ud}{\ud t}\mathcal F(g_t)
 =\{\mathcal F,\mathcal H\}(g_t)
 +\int W_t\cdot\nabla_v\mathcal F_{g_t}\ud\mu.
\]
Here the velocity differential $\nabla_v\mathcal F_g$ represents the derivative of $\mathcal F$ along perturbations of the velocity coordinates: if $u:\RR^{2d}\to\RR^d$ and $g_\varepsilon\mu=(x,v+\varepsilon u(x,v))_\#(g\mu)$, then $\frac{\ud}{\ud\varepsilon}\mathcal F(g_\varepsilon)|_{\varepsilon=0}=\int u\cdot\nabla_v\mathcal F_g\,g\ud\mu$. Such perturbations leave the spatial marginal $q\mux$ of $g\mu$, and hence the optimal transportation map $T_{q,\sigma}$, unchanged. We hence have 
\begin{equation}\label{eq:vertical-variation}
 \Cs(g_\varepsilon)-\Cs(g)
         =\varepsilon\int u\cdot\xi_{q,\sigma}\,g\ud\mu,
\end{equation}
and that the velocity differential of $\Cs$ is the Brenier displacement $\xi_{q,\sigma}$. The first-variation formulas \eqref{eq:controlled-observable} and \eqref{eq:vertical-variation}, together with the score identity \eqref{eq:score-mean}, give
\begin{equation}\label{eq:observable-structure}
 \{\Ds,\mathcal H\}=\Cs,\qquad \Gamma_v(\Ds,\mathcal F)=0,\qquad
 \Gamma_v(\Cs,\Cs)=2\Ds,\qquad \Gamma_v(\Cs,\mathcal H)=\Cs.
\end{equation}
Indeed, the first identity is \eqref{eq:controlled-observable} along Hamiltonian transport, and the second one holds since velocity variations leave the spatial marginal unchanged. The third identity is the squared $L^2(g\mu)$ norm of the velocity differential $\xi_{q,\sigma}$ of $\Cs$, and the last one follows from \eqref{eq:score-mean}:
\[
 \Gamma_v(\Cs,\mathcal H)=
 \int\xi_{q,\sigma}\cdot\nabla_v\log g\,g\ud\mu
 =\int q\xi_{q,\sigma}\cdot m\ud\mux=\Cs(g).
\]
Hence, through Cauchy--Schwarz, $\Gamma_v(\Cs,\Cs)=2\Ds$ accounts for the product of the spatial Wasserstein distance and the square root of the action in the space-time HWI inequality, and $\Gamma_v(\Cs,\mathcal H)=\Cs$ gives the damping term $\frac{\ud}{\ud t} \mathcal{D}_\sigma$ along Langevin evolution.
\end{remark}

\section{Finite-action paths and proofs of the main theorems}\label{sec:weak}

We now establish the comparison inequality \eqref{eq:smooth-comparison} rigorously under the finite-action assumptions. The proof does not use the exact identity or the regularity assumptions of Section~\ref{sec:identity}. 
Instead, the three integrals in \eqref{eq:exact-error} are handled separately: the first by the Wasserstein acceleration inequality of~\cite{LiLu2026}*{Lemma 2.7}, the second by integration by parts against the distributional Jacobian $DT_{q,\sigma}$ followed by $r-1-\log r\ge0$, and the third by the Gaussian entropy decomposition \eqref{eq:gaussian-exact}.

\subsection{The finite-action comparison inequality}

\begin{theorem}[Finite-action comparison inequality]\label{thm:weak}
Assume \eqref{eq:intro-potential}. Let $(g,W)$ satisfy \eqref{eq:controlled-intro} and \eqref{eq:intro-path}, with $g_t$ a probability density relative to $\mu$ at every time, $t\mapsto g_t\mu$ narrowly continuous on $[0,T]$, and $g_0\mu\in\PPtwo$. Then $t\mapsto g_t\mu\in C([0,T];\PPtwo)$. The entropy is finite at every time, and $t\mapsto\Ent_\mu(g_t)$ is absolutely continuous and satisfies the entropy chain rule \eqref{eq:entropy-chain} for every $0\le s\le t\le T$. For every $\sigma\in\PPtwo(\RR^d)$ of finite entropy relative to $\mux$, the function $t\mapsto\Ds(g_t)$ is continuously differentiable on $[0,T]$, with derivative $\Cs(g_t)$ of bounded variation. Moreover, it holds that in $\mathcal D'(0,T)$,
\begin{equation}\label{eq:weak-comparison}
 \frac{\ud^2}{\ud t^2}\Ds(g_t)
 +\Ent_\mu(g_t)-\Ent(\sigma\mid\mux)+K\Ds(g_t)
       \le3\Entv(g_t)+\int W_t\cdot\xi_{q_t,\sigma}\ud\mu.
\end{equation}
\end{theorem}

Throughout, compact support in phase space is understood to be uniform on each fixed closed time interval. We first treat in \Cref{lem:regular-comparison} the densities $g$ that are smooth compact perturbations of equilibrium for a potential with bounded Hessian. Their positive lower bound ensures that the conditional mean velocity and covariance are smooth, while the bounded Hessian gives at most linear growth of $\nabla U$. Outside a compact set, these densities are constant, so the corresponding laws $g\mu$ retain the tail behavior of $\mu$.
A change to a Gaussian spatial reference then yields \cref{lem:compact-comparison} for compactly supported paths without assuming a bounded Hessian, since the force correction and the change in relative entropy, together with $\Ds(g_t)$, give $\mathcal{R}_U(q_t,\sigma)$ exactly.
To complete the proof of \cref{thm:weak}, we apply this lemma to the density-current pairs constructed in Appendix~\ref{app:recovery}, which converge in the action and in the total and spatial marginal entropies. The entropy chain rule is proved independently in \cref{prop:entropy-chain}.  

\begin{lemma}[Comparison for smooth compact perturbations]\label{lem:regular-comparison}
Assume $\|\nabla^2U\|_\infty<\infty$ and $\nabla^2U\succeq K\Id$ for some $K\in\RR$. Suppose $g=\delta+f$ with $\delta>0$ constant, $f\ge0$ smooth in $(t,x,v)$ on a neighborhood of $[0,T]$ and compactly supported in phase space, and $W$ smooth and compactly supported. Assume $\int g_t\ud\mu=1$ and \eqref{eq:controlled-intro}. Then \eqref{eq:weak-comparison} holds for every $\sigma\in\PPtwo(\RR^d)$ with finite entropy relative to $\mux$.
\end{lemma}

\begin{proof}
All differential inequalities below are understood in $\mathcal D'(0,T)$, and we suppress the time variable in spatial integrals. Since $q\ge\delta$, the conditional mean velocity $m=j/q$ and its Eulerian acceleration $\partial_tm+(m\cdot\nabla_x)m$ are smooth and compactly supported. The Wasserstein acceleration inequality of~\cite{LiLu2026}*{Lemma 2.7} therefore applies, whose proof uses only the finite second moment of the fixed target. We hence have $\frac{\ud}{\ud t}\Ds(g_t)=\Cs(g_t)$ and
\begin{equation}\label{eq:weak-acceleration}
 \frac{\ud^2}{\ud t^2}\Ds(g_t)\le\int q|m|^2\ud\mux+
   \int\left[\partial_tj-\nabla_x^*\!\left(\frac{j\otimes j}{q}\right)\right]
                       \cdot\xi_{q,\sigma}\ud\mux.
\end{equation}
The inequality \eqref{eq:weak-acceleration} follows by bounding $\Ds(g_{t+s})$, for small $s$, by the cost of the coupling obtained by moving the source points of the optimal coupling at time $t$ along the flow of $m$, while keeping the target points $T_{q_t,\sigma}(x)$ fixed. This cost equals $\Ds(g_t)$ at $s=0$.

We next estimate the covariance stress using the BV integration-by-parts formula of~\cite{Lu2026}*{(5.14)}. Convexity of the Brenier potential implies that $DT_{q,\sigma}$ is a positive semidefinite matrix-valued measure, which we pair with the smooth spatial density $q(x)\e^{-U(x)}/Z_x$. In doing so, we introduce spatial cutoffs and then pass to the limit as their
supports expand. For this, we first prove
\begin{equation}\label{eq:weighted-trace}
 \int \frac{q(x)\e^{-U(x)}}{Z_x}\,
       \tr(DT_{q,\sigma})(\ud x)<\infty.
\end{equation}
Choose smooth cutoffs $0\le\chi_R\le1$ increasing to one, equal to one on $B_R$ and zero outside $B_{2R}$, with $|\nabla\chi_R|\le C/R$. Then
\[
 \int\chi_R q\mux\,\tr(DT_{q,\sigma})
       =-\int T_{q,\sigma}\cdot\nabla(\chi_R q\mux)\ud x.
\]
Since $\nabla(q\mux)=(\nabla\log q-\nabla U)\,q\mux$, the Cauchy--Schwarz inequality and $(T_{q,\sigma})_\#(q\mux)=\sigma$ give
\[
 \begin{aligned}
 &\left|\int\chi_R qT_{q,\sigma}\cdot(\nabla\log q-\nabla U)\ud\mux\right| \le\left(\int|y|^2\ud\sigma(y)\right)^{1/2}
       \left(\int q|\nabla\log q-\nabla U|^2\ud\mux\right)^{1/2},\\
 &\left|\int qT_{q,\sigma}\cdot\nabla\chi_R\ud\mux\right|
 \le\frac C R\left(\int|y|^2\ud\sigma(y)\right)^{1/2}\longrightarrow0.
 \end{aligned}
\]
The integral $\int q|\nabla\log q-\nabla U|^2\ud\mux$ is finite: outside a compact set, $q$ is constant, and the bounded Hessian implies that $\nabla U$ grows at most linearly. On compact sets, $q$ is smooth and bounded below by $\delta$. Then, the monotone convergence theorem proves \eqref{eq:weighted-trace}. 

For $g=\delta+f$, the covariance $\Sigma_g$ is smooth, bounded, and uniformly positive definite. Applying the BV integration-by-parts formula \cite{Lu2026}*{(5.14)} to $\chi_R q\Sigma_g$ and passing to the limit using \eqref{eq:weighted-trace} yields
\begin{equation}\label{eq:bv-stress}
 \int\xi_{q,\sigma}\cdot\nabla_x^*(q\Sigma_g)\ud\mux  = \int q\tr\Sigma_g\ud\mux - \int\frac{q(x)\e^{-U(x)}}{Z_x}\,    \Sigma_g(x):DT_{q,\sigma}(\ud x).
\end{equation}
Indeed, the cutoff error is bounded by $C R^{-1}\int_{R<|x|<2R}q|\xi_{q,\sigma}|\ud\mux$ and tends to zero. Let $G$ denote the density of the absolutely continuous part of $DT_{q,\sigma}$. Since the singular part is positive semidefinite, its contribution to \eqref{eq:bv-stress} is nonpositive. Hence, 
\begin{equation} \label{eq:bv-stress2}
     \int\xi_{q,\sigma}\cdot\nabla_x^*(q\Sigma_g)\ud\mux \le \int q\bigl[\tr\Sigma_g-\tr(\Sigma_gG)\bigr]\ud\mux.     
\end{equation}

It remains to identify the entropy terms. Since both $q\mux$ and $\sigma$ are absolutely continuous, the almost-everywhere transport Jacobian formula~\cite{villani2021topics}*{Theorem 4.8} gives $\det G>0$ and $q\mux=\sigma(T_{q,\sigma})\det G$ for $q\mux$-almost every $x$, as in the proof of \cref{thm:exact}.

The functions $\log q$ and $\log\frac{\sigma}{\mux}(T_{q,\sigma})$ are integrable against $q\mux$ because $\Ent_{\mux}(q)$ and $\Ent(\sigma\mid\mux)$ are finite. The bounded Hessian gives $|U(x)|\le C(1+|x|^2)$. Since $q\mux$ and $\sigma$ have finite second moments, $U$ and $U(T_{q,\sigma})$ are also integrable against $q\mux$. Hence $\log\det G$ is integrable against $q\mux$. Dividing the Jacobian relation by $\mux=\e^{-U}/Z_x$ and integrating its logarithm yields \eqref{eq:transport-entropy}.

The weighted trace bound \eqref{eq:weighted-trace} also ensures integrability of $\tr(\Sigma_gG)$, while $\log\det\Sigma_g$ is bounded. The conditional mean and covariance terms satisfy
\[
 \begin{aligned}
 &|m|^2+\tr\Sigma_g-\tr(\Sigma_gG)+\log\det G\\
 &\qquad\le |m|^2+\tr\Sigma_g-d-\log\det\Sigma_g
 \le2\Ent_\kv(h_x).
 \end{aligned}
\]
The first inequality follows by applying $r-1-\log r\ge0$ to the eigenvalues of $\Sigma_g^{1/2}G\Sigma_g^{1/2}$, and the second follows from \eqref{eq:gaussian-exact}. Substituting \eqref{eq:moment-equations} and \eqref{eq:bv-stress2} into \eqref{eq:weak-acceleration}, and using \eqref{eq:transport-entropy}, we obtain
\[
 \frac{\ud^2}{\ud t^2}\Ds(g_t)+\Ent_{\mux}(q)-\Ent(\sigma\mid\mux)+\mathcal R_U(q,\sigma)
          \le2\Entv(g)+\int W\cdot\xi_{q,\sigma}\ud\mu.
\]
Finally, we add $\Entv(g)$ to both sides, use the entropy decomposition \eqref{eq:conditional-entropy}, and apply $\mathcal R_U\ge K\Ds$ to obtain the desired inequality \eqref{eq:weak-comparison}.
\end{proof}

\begin{lemma}[Comparison for compactly supported paths]\label{lem:compact-comparison}
Assume \eqref{eq:intro-potential}. Let $(\varrho,F)$ be a smooth solution of \eqref{eq:ordinary-control}, of mass one, with density and current compactly supported uniformly on a neighborhood of $[0,T]$, and with finite quadratic action. Then \eqref{eq:weak-comparison} holds for every $\sigma\in\PPtwo(\RR^d)$ with $\Ent(\sigma\mid\mux)<\infty$.
\end{lemma}
\begin{proof}
We let $\widetilde U(x)=|x|^2/2$ and $\widetilde\mu_x$ be its Gaussian Gibbs law, and we write $\widetilde\mu=\widetilde\mu_x\otimes\kv$ and $\widetilde\varrho_\infty$ for its Lebesgue density. 
We will prove the comparison inequality relative to $\widetilde\mu$, where \cref{lem:regular-comparison} applies, and then return to the original reference measure $\mu$. 
To transfer the controlled equation from $U$ to $\widetilde U$, set
\[
 \widetilde F=F+(\nabla\widetilde U-\nabla U)\varrho.
\]
Then $(\varrho,\widetilde F)$ solves \eqref{eq:ordinary-control} with $\widetilde U$ in place of $U$. Since $\nabla\widetilde U-\nabla U$ is bounded on the compact source support, $\widetilde F/\sqrt\varrho\in L^2(\ud t\ud z)$. Applying \cref{lem:energy-integrability} to $\sigma\otimes\kv$ gives $U\in L^1(\sigma)$. Since $\widetilde U$ is quadratic and $\sigma\in\PPtwo$, the change of reference in relative entropy gives $\Ent(\sigma\mid\widetilde\mu_x)<\infty$.

To apply \cref{lem:regular-comparison}, we mix the density with equilibrium: for $\delta\in(0,1)$, 
\[
 \varrho_\delta=(1-\delta)\varrho+\delta\widetilde\varrho_\infty,\qquad
 F_\delta=(1-\delta)\widetilde F.
\]
The pair still solves the controlled equation for $\widetilde U$. Relative to $\widetilde\mu$, its density is a positive constant plus a smooth nonnegative compact perturbation, and its velocity current is smooth and compactly supported. Thus \cref{lem:regular-comparison} gives \eqref{eq:weak-comparison} with $\widetilde U$ in place of $U$, $K=1$, and target $\sigma$.

Let $\delta\downarrow0$. The laws converge uniformly in $W_2$, and their densities converge in space-time $L^1$. The total entropies converge in $L^1(0,T)$ by the mixture bounds
\[
 (1-\delta)\Ent(\nu\mid\widetilde\mu)+\delta\log\delta+(1-\delta)\log(1-\delta)
 \le\Ent((1-\delta)\nu+\delta\widetilde\mu\mid\widetilde\mu)
 \le(1-\delta)\Ent(\nu\mid\widetilde\mu),
\]
applied to $\nu=\varrho_t\ud z$. Applying the same bounds to the spatial marginals gives convergence of the marginal entropies, and hence of the conditional velocity entropies. Moreover, $F_\delta/\sqrt{\varrho_\delta}\to\widetilde F/\sqrt\varrho$ strongly in $L^2$ by dominated convergence. Together with \cref{lem:transport-stability}, this proves convergence of the control pairing. Passing to the limit in distributions gives the comparison for $(\varrho,\widetilde F)$ relative to $\widetilde\mu$.

We now return to the original reference measure. The spatial law, its optimal transportation map, and the conditional velocity entropy are fixed. Write $\nu_t=\varrho_t\ud z=g_t\mu$. The normalizing constants cancel in the entropy difference, and the force correction combines with the change of relative entropy to give
\[
 \begin{aligned}
 &\Ent(\nu_t\mid\widetilde\mu)-\Ent(\sigma\mid\widetilde\mu_x)
       +\Ds(g_t)-\int\widetilde F_t\cdot\xi_{q_t,\sigma}\ud z\\
 &\quad=\Ent_\mu(g_t)-\Ent(\sigma\mid\mux)
       +\mathcal R_U(q_t,\sigma)-\int F_t\cdot\xi_{q_t,\sigma}\ud z.
 \end{aligned}
\]
Here we used $\widetilde U(y)-\widetilde U(x)-\nabla\widetilde U(x)\cdot(y-x)=|x-y|^2/2$. All terms are integrable: the source has compact support, while $U\in L^1(\sigma)$ and $\sigma$ has a finite second moment. Rewriting the comparison for $\widetilde U$ in terms of $\mu$ and $F$, and using $\mathcal R_U\ge K\Ds$, proves \eqref{eq:weak-comparison}.
\end{proof}

\begin{proof}[Proof of \cref{thm:weak}]
Lemma~\ref{lem:superposition}, applied with $W=gA$ as in \eqref{eq:ordinary-control}, first gives $t\mapsto\nu_t=g_t\mu\in C([0,T];\PPtwo)$. Proposition~\ref{prop:entropy-chain} gives finiteness and absolute continuity of the entropy, together with \eqref{eq:entropy-chain}. We now identify the derivative of $\Ds(g_t)$. The spatial continuity equation in \eqref{eq:moment-equations} has velocity $m=j/q$, and
\[
 \int_0^T\!\int q_t|m_t|^2\ud\mux\ud t
 \le\int_0^T\!\int|v|^2\ud\nu_t\ud t<\infty,
\]
where finiteness follows from the $W_2$ continuity of $\nu_t=g_t\mu$ on $[0,T]$. The Wasserstein first-variation formula \cite{ambrosio2005gradient}*{Theorem 8.4.7} therefore gives $\frac{\ud}{\ud t}\Ds(g_t)=\Cs(g_t)$ almost everywhere. By \cref{lem:transport-stability}, $\Cs(g_t)$ is continuous on $[0,T]$. Hence $\Ds\in C^1([0,T])$, and the identity holds at every time.

We next prove \eqref{eq:weak-comparison} for the given comparison measure $\sigma$. Let $(\varrho^k,F^k)$ be the approximation from \cref{thm:recovery}, and write $g^k=\varrho^k/\varrho_\infty$ and $q^k=\Piv g^k$. Lemma~\ref{lem:compact-comparison} gives the comparison inequality for each approximating pair. The total and spatial marginal entropies converge in $L^1(0,T)$, so the conditional velocity entropies also converge in $L^1(0,T)$. The transport distances converge uniformly on $[0,T]$. For the control term, \eqref{eq:displacement-strong} and strong convergence of $F^k/\sqrt{\varrho^k}$ give
\[
 F^k\cdot\xi_{q^k,\sigma}
 =\frac{F^k}{\sqrt{\varrho^k}}\cdot
   \bigl(\sqrt{\varrho^k}\,\xi_{q^k,\sigma}\bigr)
 \longrightarrow F\cdot\xi_{q,\sigma}
 \quad\text{in }L^1(\ud t\ud z).
\]
Testing the comparison inequality with a nonnegative smooth function compactly supported in $(0,T)$ and passing to the limit in these terms proves \eqref{eq:weak-comparison} for $\sigma$.

It remains to prove bounded variation of $\Cs(g_t)$ and finiteness of the remainder measure. Every term in \eqref{eq:weak-comparison} other than $\frac{\ud}{\ud t}\Cs(g_t)$ belongs to $L^1(0,T)$. In particular,
\[
 \left|\int W_t\cdot\xi_{q_t,\sigma}\ud\mu\right|
 \le W_2(q_t\mux,\sigma)
 \left(\int\frac{|W_t|^2}{g_t}\ud\mu\right)^{1/2}.
\]
The distance factor is bounded on $[0,T]$, and the action is finite, so this bound is integrable in time. Thus \eqref{eq:weak-comparison} has the form $\frac{\ud}{\ud t}\Cs(g_t)\le f(t)$ for some $f\in L^1(0,T)$. The continuous function $t\mapsto\Cs(g_t)-\int_0^t f(r)\ud r$ is nonincreasing and has finite endpoint values. It therefore has finite total variation, which proves $\Cs(g_\cdot)\in BV([0,T])$ and finiteness of the nonnegative remainder measure.
\end{proof}

\begin{remark}[Sharpness of the comparison inequality \eqref{eq:weak-comparison}]\label{rem:sharp-three}
We take $U(x)=x^2/2$ in one dimension and define the path of densities: 
\[
 g_t\mu=\mathcal N(t-1,1)\otimes\mathcal N(1,1),
 \qquad 0\le t\le2.
\]
Direct differentiation gives
\[
 (\partial_t+\La)g_t=(t-1)(v-1)g_t=\nabla_v^*W_t\,,  \qquad W_t=(t-1)g_t\,.
\]
This smooth path satisfies the assumptions of \cref{thm:weak}. For comparison with $\sigma=\mux$, at $t=1$ we have
\[
 \begin{gathered}
 W_1=0,\qquad \Ds(g_1)=0,\qquad
 \frac{\ud^2}{\ud t^2}\Ds(g_t)\big|_{t=1}=1,\\
 \Ent_\mu(g_1)=\Entv(g_1)=1/2.
 \end{gathered}
\]
All terms in \eqref{eq:weak-comparison} are continuous along this path, so the inequality holds at every time. At $t=1$, the terms $\Ent(\sigma\mid\mux)$, $K\Ds(g_1)$, and $\int W_1\cdot\xi_{q_1,\sigma}\ud\mu$ vanish, and \eqref{eq:weak-comparison} with a coefficient $c$ in place of $3$ reads
\[
 \frac{\ud^2}{\ud t^2}\Ds(g_t)\big|_{t=1}+\Ent_\mu(g_1)=\frac32
 \le c\,\Entv(g_1)=\frac c2.
\]
Hence $c\ge3$. With $K=1$ and $c=3$, equality holds at every time, since the path satisfies all the conditions in \cref{rem:nonnegative-terms}: the optimal transportation map and its derivatives satisfy
\[
 \begin{gathered}
 T_t(x)=x-(t-1),\qquad G_t=\Id=\Sigma_{g_t}^{-1},\\
 \partial_tT_t+G_tm_t=-1+1=0,
 \end{gathered}
\]
the conditional velocity laws are Gaussian, and $\mathcal R_U=\Ds$ for the quadratic potential $U$.
\end{remark}

\subsection{Time integration and the space-time HWI inequality}\label{subsec:hwi}
As a consequence of the finite-action comparison inequality in \cref{thm:weak}, we next establish the space-time HWI inequality for a general spatial comparison measure, extending \cref{thm:main-hwi}.  For a fixed comparison measure $\sigma$, write
\[
 W_{x,\sigma,T}^2(g)=\frac1T\int_0^T W_2^2(q_t\mux,\sigma)\ud t.
\]
Since $t\mapsto\Cs(g_t)$ is continuous on $[0,T]$, its endpoint traces coincide with $\Cs(g_0)$ and $\Cs(g_T)$.

\begin{theorem}[HWI with a spatial comparison measure]\label{thm:comparator-hwi}
Let $(g,W)$ satisfy the assumptions of \cref{thm:weak}, and fix $\sigma\in\PPtwo(\RR^d)$ with finite entropy relative to $\mux$. Then
\begin{align}\label{eq:comparator-hwi}
 \Ent_{\muT}(g)-\Ent(\sigma\mid\mux)\le{}&
 \frac{\Cs(g_0)-\Cs(g_T)}{T}+W_{x,\sigma,T}(g)\sqrt{\mathcal A_T(g)}
                +3W_{v,T}(g)\sqrt{I_{v,T}(g)}\nonumber\\
 &-\frac K2W_{x,\sigma,T}^2(g)-\frac32W_{v,T}^2(g).
\end{align}
\end{theorem}

\begin{proof}
By \cref{thm:weak}, $\frac{\ud}{\ud t}\Ds(g_t)=\Cs(g_t)$ and $\Cs(g_t)$ is continuous and of bounded variation. We therefore integrate \eqref{eq:weak-comparison} over $(0,T)$; after division by $T$, the derivative contributes the boundary term $(\Cs(g_0)-\Cs(g_T))/T$. The control term is bounded by 
\[
 \frac1T\int_0^T\!\int W_t\cdot\xi_{q_t,\sigma}\ud\mu\ud t
       \le W_{x,\sigma,T}(g)\sqrt{J_T(g,W)}.
\]
For the conditional velocity entropy, using Gaussian HWI \eqref{eq:fiber-hwi} followed by the Cauchy--Schwarz inequality in time, we have 
\[
 \frac1T\int_0^T\Entv(g_t)\ud t
       \le W_{v,T}(g)\sqrt{I_{v,T}(g)}-\tfrac12W_{v,T}^2(g).
\]
Substituting these bounds into the integrated comparison inequality gives \eqref{eq:comparator-hwi} with $J_T(g,W)$ in place of $\mathcal A_T(g)$. For the fixed path $g$, \cref{thm:weak} applies to every driving current of finite quadratic action. Taking the infimum over these currents then proves \eqref{eq:comparator-hwi}.
\end{proof}

In \eqref{eq:comparator-hwi}, the comparison measure enters through its entropy relative to $\mux$ and the optimal transportation map to $\sigma$. Thus, the left-hand side is a difference of entropies relative to the same reference measure $\mu$, rather than an entropy relative to $\sigma\otimes\kv$.

\begin{remark}[Static restriction]\label{rem:static}
Let $g_t(x,v)=q(x)$ be independent of $t$ and $v$, with $q\mux\in\PPtwo(\RR^d)$, finite entropy, and finite spatial Fisher information $I_x(q):=\int|\nabla_xq|^2/q\ud\mux$. Then $W=\nabla_xq$ drives the path, since
$\La q=v\cdot\nabla_xq=\nabla_v^*(\nabla_xq)$, and the pair has quadratic action $J_T(g,W)=I_x(q)$. Conversely, for any velocity current $W$ of finite quadratic action driving this path, taking the first velocity moment of \eqref{eq:controlled-intro}, justified by velocity cutoffs, gives $\Piv W_t=\nabla_xq$ for almost every $t$. It follows that
\[
 \int\frac{|W_t|^2}{q}\ud\mu
    \ge\int\frac{|\Piv W_t|^2}{q}\ud\mux
    =I_x(q).
\]
Averaging in time and minimizing over $W$ proves $\mathcal A_T(g)=I_x(q)$. Moreover, note $W_{v,T}=I_{v,T}=0$, and $\Cs(g_0)=\Cs(g_T)$. Thus \eqref{eq:comparator-hwi} restricts to the usual spatial HWI (and, for $\sigma=\mux$, to \eqref{eq:classical-hwi}) wherever the assumptions of \cref{thm:comparator-hwi} hold. 
\end{remark}

\begin{remark}[Necessity of the boundary term]\label{rem:endpoint}
The boundary term in \cref{thm:main-hwi} cannot be omitted. In fact, Hamiltonian transport preserves entropy, while a path that initially differs from equilibrium only in position can have small velocity Fisher information on short intervals. 

To see this, take $U(x)=Kx^2/2$ in one dimension, with $K>0$, fix
$r>0$, and define
\begin{gather*}
 g_t\mu=\mathcal N(X_t,K^{-1})\otimes\mathcal N(V_t,1),\\
 X_t=r\cos(\sqrt K\,t),\qquad
 V_t=-r\sqrt K\sin(\sqrt K\,t).
\end{gather*}
Since $\frac{\ud X_t}{\ud t}=V_t$ and $\frac{\ud V_t}{\ud t}=-KX_t$, this path satisfies \eqref{eq:controlled-intro} with $W=0$. Thus $\mathcal A_T(g)=0$, and 
\begin{equation} \label{rem:axueq}
    \Ent_\mu(g_t)=\frac12\bigl(KX_t^2+V_t^2\bigr) = \frac{Kr^2}{2}.
\end{equation}
The initial conditional velocity laws are at equilibrium, and $\Iv(g_t)=V_t^2$ tends to zero as $t\downarrow0$. Consequently, as $T\downarrow0$, we have $W_{x,T}^2(g)\to r^2$ and $W_{v,T}^2(g)=I_{v,T}(g)\to 0$. Without the boundary term, the right-hand side of \eqref{eq:main-hwi} would therefore tend to $-Kr^2/2$, whereas the left-hand side remains equal to $Kr^2/2$ by \eqref{rem:axueq}. 

In fact, \eqref{eq:main-hwi} holds with equality along this path for every $T>0$. Setting $s_T=\sin(2\sqrt K\,T)/(2\sqrt K\,T)$, a direct computation gives
\[
 \frac{\mathcal C_{\mux}(g_0)-\mathcal C_{\mux}(g_T)}{T}  =Kr^2s_T,
 \qquad W_{x,T}^2(g)=\frac{r^2}{2}(1+s_T),
 \qquad W_{v,T}^2(g)=I_{v,T}(g)=\frac{Kr^2}{2}(1-s_T).
\]
Substitution into \eqref{eq:main-hwi} gives $Kr^2/2$ on both sides. 
In addition, note that the boundary-free LSI in \cref{thm:long-lsi} holds on every positive time interval, with a coefficient of the velocity Fisher information growing as $T^{-2}$ at short times. The same path shows that this growth is necessary; see \cref{rem:lsi-short-time}. 
\end{remark}

\subsection{The kinetic EVI and entropy balance}\label{subsec:langevin-balance}

We next specialize the comparison inequality to Langevin evolution to prove \cref{thm:kinetic-evi}. The velocity Ornstein--Uhlenbeck generator $\gamma \mathcal{L}_s$ produces the damping term, and the entropy chain rule becomes the entropy balance. The following proposition, proved in Appendix~\ref{app:langevin-paths}, verifies the assumptions needed to apply \cref{thm:weak} to a Langevin path starting from a density of finite relative entropy and finite second moment.

\begin{proposition}[Langevin paths and entropy balance]\label{prop:semigroup-paths}
Assume \eqref{eq:intro-potential} and fix $\gamma>0$. The Langevin equation
\[
 \ud X_t=V_t\ud t,\qquad
 \ud V_t=-\nabla U(X_t)\ud t-\gamma V_t\ud t+\sqrt{2\gamma}\ud B_t
\]
is nonexplosive and preserves $\mu$. If $\nu_0=g\mu\in\PPtwo$ and $\Ent_\mu(g)<\infty$, its law $\nu_t=(\Pt_tg)\mu$ belongs to $C([0,T];\PPtwo)$ for every finite $T$, has a density at every time, and satisfies
\begin{equation}\label{eq:semigroup-entropy-balance}
 \Ent_\mu(\Pt_tg)+\gamma\int_0^t \Iv(\Pt_rg)\ud r=\Ent_\mu(g).
\end{equation}
\end{proposition}

We now give the proof of Kinetic EVI in \cref{thm:kinetic-evi}. 

\begin{proof}[Proof of \cref{thm:kinetic-evi}]
Fix $T>0$ and set $g_t=\Pt_tg$. By \cref{prop:semigroup-paths}, the path $g_t\mu$ is continuous in $W_2$ and satisfies the entropy balance \eqref{eq:semigroup-entropy-balance}. Its velocity current is
\begin{equation}\label{eq:physical-control}
 W_t=-\gamma\nabla_vg_t,
 \qquad
 \int_0^T\!\int\frac{|W_t|^2}{g_t}\ud\mu\ud t
                         =\gamma^2\int_0^T \Iv(g_t)\ud t<\infty.
\end{equation}
Thus the Langevin path satisfies the assumptions of \cref{thm:weak} on $[0,T]$. For almost every $t$, the conditional mean identity \eqref{eq:score-mean} gives
\[
 \int W_t\cdot\xi_{q_t,\sigma}\ud\mu
 =-\gamma\Cs(g_t)=-\gamma\frac{\ud}{\ud t}\Ds(g_t).
\]
Substituting this pairing into \eqref{eq:weak-comparison} gives \eqref{eq:intro-langevin-comparison} in $\mathcal D'(0,T)$, with the damping term $\gamma\frac{\ud}{\ud t}\Ds(g_t)$, where the second inequality is by the conditional Gaussian LSI. The regularity of $\Ds(g_t)$ and $\Cs(g_t)$ follows from \cref{thm:weak}. Since $T>0$ was arbitrary, the conclusions hold on $[0,\infty)$.
\end{proof}

\section{Space-time logarithmic Sobolev inequalities}\label{sec:inequalities}

We derive space-time LSI by controlling the transport distances in the HWI and comparison inequalities. We first assume $\mathcal C_{\mux}(g_T)\ge\mathcal C_{\mux}(g_0)$, in which case the boundary term in \eqref{eq:main-hwi} is nonpositive, and we see that the coefficients in the space-time LSI do not depend on the length of the time interval (\cref{cor:boundary-lsi}). 
We then remove this restriction on every positive time interval; without it, the coefficient of velocity Fisher information must grow at short times, by \cref{rem:lsi-short-time}.

\subsection{Inequalities with a nonpositive boundary term}

\begin{corollary}[Space-time LSI with a nonpositive boundary term]\label{cor:boundary-lsi}
Under the assumptions of \cref{thm:comparator-hwi}, with lower Hessian bound $\nabla^2U\succeq K\Id$, comparison measure $\sigma=\mux$, and $\mathcal C_{\mux}(g_T)\ge\mathcal C_{\mux}(g_0)$, the following estimates hold.
If $K>0$, then
\begin{equation}\label{eq:positive-lsi}
 \Ent_{\muT}(g)\le\frac{1}{2K}\mathcal A_T(g)+\frac32I_{v,T}(g).
\end{equation}
If $K\le0$, $\mux$ satisfies Talagrand's transportation inequality \eqref{eq:spatial-t2}, and $\alpha:=1+K/\rho>0$, then 
\begin{align}\label{eq:negative-lsi}
 \Ent_{\muT}(g)
 &\le\frac12I_{v,T}(g)
 +\frac{\bigl[\sqrt{\mathcal A_T(g)}
       +\sqrt{\mathcal A_T(g)+2\alpha\rho I_{v,T}(g)}\bigr]^2}
             {2\alpha^2\rho}\nonumber\\
 &\le\frac{2}{\alpha^2\rho}\mathcal A_T(g)
                   +\left(\frac12+\frac2\alpha\right)I_{v,T}(g).
\end{align}
\end{corollary}

\begin{proof}
Suppose first that $K>0$. We discard the nonpositive boundary term in \eqref{eq:comparator-hwi} and use
\[
 W_{x,T}(g)\sqrt{\mathcal A_T(g)}-\frac K2W_{x,T}^2(g)\le\frac{\mathcal A_T(g)}{2K},\qquad
 3W_{v,T}(g)\sqrt{I_{v,T}(g)}-\frac32W_{v,T}^2(g)\le\frac32I_{v,T}(g),
\]
by completing the square. Then \eqref{eq:positive-lsi} follows.

For $K\le0$, we return to the comparison inequality \eqref{eq:weak-comparison} with $\sigma=\mux$. Integrating in time, using the entropy decomposition \eqref{eq:conditional-entropy}, and discarding the nonpositive boundary term gives, after applying the Cauchy--Schwarz inequality and minimization over $W$, 
\begin{equation} \label{eq:auxeqlsi}
     \frac1T\int_0^T\Ent_{\mux}(q_t)\ud t  +\frac K2W_{x,T}^2(g)
       \le\frac2T\int_0^T\Entv(g_t)\ud t +W_{x,T}(g)\sqrt{\mathcal A_T(g)}.
\end{equation}
Averaging the transportation inequality \eqref{eq:spatial-t2} in time gives
\begin{equation} \label{eq:auxeqlsi2}
     W_{x,T}^2(g)\le\frac2{\rho T}\int_0^T\Ent_{\mux}(q_t)\ud t.
\end{equation}
Since $K\le0$, this bound gives, for the left-hand side of \eqref{eq:auxeqlsi},
\[
 \frac1T\int_0^T\Ent_{\mux}(q_t)\ud t+\frac K2W_{x,T}^2(g)
 \ge\Bigl(1+\frac K\rho\Bigr)\frac1T\int_0^T\Ent_{\mux}(q_t)\ud t =\frac\alpha T\int_0^T\Ent_{\mux}(q_t)\ud t.
\]
On the right-hand side, the same bound \eqref{eq:auxeqlsi2} controls $W_{x,T}(g)$, and the Gaussian LSI $\Entv(g_t)\le\frac12\Iv(g_t)$ gives $\frac2T\int_0^T\Entv(g_t)\ud t\le I_{v,T}(g)$. Hence, we have 
\[
 \frac\alpha T\int_0^T\Ent_{\mux}(q_t)\ud t
 \le I_{v,T}(g)
 +\left(\frac{2\mathcal A_T(g)}{\rho T}\int_0^T\Ent_{\mux}(q_t)\ud t\right)^{1/2}.
\]
Solving this quadratic inequality for the square root of the time-averaged spatial entropy yields
\[
 \frac1T\int_0^T\Ent_{\mux}(q_t)\ud t
 \le\frac{\bigl[\sqrt{\mathcal A_T(g)}
       +\sqrt{\mathcal A_T(g)+2\alpha\rho I_{v,T}(g)}\bigr]^2}
      {2\alpha^2\rho}.
\]
The entropy decomposition and Gaussian LSI give $\Ent_{\muT}(g)\le T^{-1}\int_0^T\Ent_{\mux}(q_t)\ud t+I_{v,T}(g)/2$, proving the first inequality in \eqref{eq:negative-lsi}. The second follows from $(a+b)^2\le2a^2+2b^2$, completing the proof. 
\end{proof}

\subsection{A boundary-free inequality on arbitrary intervals}

  We now prove a space-time LSI on every positive time interval, without assuming $\mathcal C_{\mux}(g_T)\ge\mathcal C_{\mux}(g_0)$. Its coefficient of velocity Fisher information remains bounded on long intervals and has the optimal inverse-square order on short intervals.

\begin{theorem}[Boundary-free space-time LSI]\label{thm:long-lsi}
Suppose that $\nabla^2U\succeq K\Id$ for some $K\le0$, that $\mux$ satisfies \eqref{eq:spatial-t2}, and that $\alpha=1+K/\rho>0$. Let $(g,W)$ satisfy the assumptions of \cref{thm:weak}. There is a universal constant $C$ such that, for every $T>0$,
\begin{equation}\label{eq:long-lsi}
 \Ent_{\muT}(g)\le\frac{C}{\alpha^2}
                \left[\max\left\{1,\frac1{\rho T^2}\right\}
                       I_{v,T}(g)+\rho^{-1}\mathcal A_T(g)\right].
\end{equation}
In particular, the coefficient of $I_{v,T}(g)$ is independent of $T$ for $T\ge\rho^{-1/2}$. For a convex potential, one may take $K=0$ and $\alpha=1$.
\end{theorem}

\begin{proof}
We fix a current $W$ driving the path $g$ with $J_T(g,W)<\infty$. We first estimate the entropy on $[0,L]\subset[0,T]$ and then apply the estimate on a partition of $[0,T]$, following the time-localization argument of~\cite{LiLu2026}*{Theorem 3.7 and Corollary 3.14}. In this proof, we write $\Ent_{\bar\mu_L}(g)$ for the entropy averaged over $[0,L]$, and abbreviate $J_L=J_L(g,W)$ and $I_{v,L}=I_{v,L}(g)$, as in \eqref{eq:intro-averages}. Below $C$ denotes a universal constant that may change from line to line.

The assumptions give the following bounds. The entropy derivative and the estimates at a fixed time hold almost everywhere, while the differential inequality for $\mathcal C_{\mux}$ holds in $\mathcal D'(0,T)$:
\begin{equation} \label{eq:scalar-controlled}
    \begin{aligned}
     & \Bigl|\frac{\ud}{\ud t}\Ent_\mu(g_t)\Bigr|
   \le\Bigl(\int\frac{|W_t|^2}{g_t}\ud\mu\Bigr)^{1/2}\Iv(g_t)^{1/2}, \\
& \frac{\ud}{\ud t}\mathcal C_{\mux}(g_t)+\alpha\Ent_\mu(g_t)
   \le\frac{2+\alpha}{2}\Iv(g_t)+\int W_t\cdot\xi_{q_t}\ud\mu, \\
& |\mathcal C_{\mux}(g_t)|
   \le\Bigl(\frac2\rho\Ent_\mu(g_t)\Iv(g_t)\Bigr)^{1/2},\qquad
 \Bigl|\int W_t\cdot\xi_{q_t}\ud\mu\Bigr|
   \le\Bigl(\frac2\rho\Ent_\mu(g_t)\int\frac{|W_t|^2}{g_t}\ud\mu\Bigr)^{1/2}.    
\end{aligned}    
\end{equation}
The first bound follows from the entropy chain rule and the Cauchy--Schwarz inequality. To prove the second, use the entropy decomposition \eqref{eq:conditional-entropy} to write the comparison inequality \eqref{eq:weak-comparison} with $\sigma = \mu_x$ as
\[
 \frac{\ud}{\ud t}\mathcal C_{\mux}(g_t)
 +\Ent_{\mux}(q_t)+\frac K2W_2^2(q_t\mux,\mux)
 \le2\Entv(g_t)+\int W_t\cdot\xi_{q_t}\ud\mu.
\]
Since $K\le0$, the transportation inequality gives
\[
 \Ent_{\mux}(q_t)+\frac K2W_2^2(q_t\mux,\mux)
 \ge\alpha\Ent_{\mux}(q_t).
\]
Adding $\alpha\Entv(g_t)$ to both sides and applying the Gaussian LSI proves the second line of \eqref{eq:scalar-controlled}. The bound on $\mathcal C_{\mux}$ is \eqref{eq:current-t2}; the bound on $\int W_t\cdot\xi_{q_t}\ud\mu$ follows from Cauchy--Schwarz and \eqref{eq:spatial-t2}.

Integrating the first line of \eqref{eq:scalar-controlled} and applying Cauchy--Schwarz in time bounds the deviation of $\Ent_\mu(g_t)$ from its time average by $L\sqrt{J_LI_{v,L}}$:
\begin{equation}\label{eq:entropy-oscillation}
 \bigl|\Ent_\mu(g_t)-\Ent_{\bar\mu_L}(g)\bigr|\le L\sqrt{J_LI_{v,L}},
 \qquad 0\le t\le L.
\end{equation}
To remove the boundary term, we fix a nonnegative $\chi\in C_c^\infty(0,1)$ with integral one and define the rescaled time cutoff
\begin{equation}\label{eq:time-cutoff}
 \chi_L(t)=\frac1L\chi\!\left(\frac tL\right),
\end{equation}
which is compactly supported in $(0,L)$, and satisfies
\begin{equation}\label{eq:time-cutoff-norms}
 \|\chi_L\|_\infty\le\frac C{L},\quad
 \int_0^L\chi_L^2\le\frac C{L},\quad
 \int_0^L\left(\frac{\ud\chi_L}{\ud t}(t)\right)^2\ud t\le\frac C{L^3}.
\end{equation}

We next test the second line of \eqref{eq:scalar-controlled} with $\chi_L$. Integration by parts gives
\[
 \alpha\int_0^L\chi_L\Ent_\mu(g_t)\ud t
 \le\frac{2+\alpha}{2}\int_0^L\chi_L\Iv(g_t)\ud t
       +\int_0^L\chi_L\!\int W_t\cdot\xi_{q_t}\ud\mu\ud t
       +\int_0^L\frac{\ud\chi_L}{\ud t}(t)\mathcal C_{\mux}(g_t)\ud t.
\]
By \eqref{eq:entropy-oscillation}, we have
\[
 \begin{aligned}
 \left|\int_0^L\chi_L\Ent_\mu(g_t)\ud t-\Ent_{\bar\mu_L}(g)\right|
 &\le L\sqrt{J_LI_{v,L}},\\
 \Ent_\mu(g_t)&\le\Ent_{\bar\mu_L}(g)+L\sqrt{J_LI_{v,L}}.
 \end{aligned}
\]
Using the Cauchy--Schwarz inequality in time, \eqref{eq:scalar-controlled}, and the cutoff bounds \eqref{eq:time-cutoff-norms}, we obtain
\[
 \begin{aligned}
 &\left|\int_0^L\chi_L\!\int W_t\cdot\xi_{q_t}\ud\mu\ud t\right| \le\frac C{\sqrt\rho}
    \bigl(\Ent_{\bar\mu_L}(g)+L\sqrt{J_LI_{v,L}}\bigr)^{1/2}\sqrt{J_L},\\
 &\left|\int_0^L\frac{\ud\chi_L}{\ud t}(t)\mathcal C_{\mux}(g_t)\ud t\right| \le\frac C{L\sqrt\rho}
    \bigl(\Ent_{\bar\mu_L}(g)+L\sqrt{J_LI_{v,L}}\bigr)^{1/2}\sqrt{I_{v,L}}.
 \end{aligned}
\]
Since $0<\alpha\le1$, these estimates yield
\[
 \alpha\Ent_{\bar\mu_L}(g)
 \le\alpha L\sqrt{J_LI_{v,L}}+CI_{v,L}
       +\frac C{\sqrt\rho}\bigl(\Ent_{\bar\mu_L}(g)+L\sqrt{J_LI_{v,L}}\bigr)^{1/2}
                       \left(\sqrt{J_L}+\frac{\sqrt{I_{v,L}}}{L}\right).
\]
Then, by Young's inequality, we have 
\[
 \begin{aligned}
 &\frac C{\sqrt\rho}\bigl(\Ent_{\bar\mu_L}(g)+L\sqrt{J_LI_{v,L}}\bigr)^{1/2}
       \left(\sqrt{J_L}+\frac{\sqrt{I_{v,L}}}{L}\right) \le\frac\alpha2\bigl(\Ent_{\bar\mu_L}(g)+L\sqrt{J_LI_{v,L}}\bigr)
       +\frac C{\alpha\rho}\left(J_L+\frac{I_{v,L}}{L^2}\right).
 \end{aligned}
\]
Absorbing $\alpha\Ent_{\bar\mu_L}(g)/2$ yields
\begin{equation}\label{eq:window-estimate}
 \Ent_{\bar\mu_L}(g)\le C\left[
 L\sqrt{J_LI_{v,L}}+\frac{I_{v,L}}\alpha
       +\frac1{\alpha^2\rho}\left(J_L+\frac{I_{v,L}}{L^2}\right)\right].
\end{equation}
For $0<L\le2\rho^{-1/2}$, note 
$L\sqrt{J_LI_{v,L}}\le C(I_{v,L}+J_L/\rho)$. Since $0<\alpha\le1$, \eqref{eq:window-estimate} yields
\[
 \Ent_{\bar\mu_L}(g)
 \le\frac{C}{\alpha^2}\left[\max\left\{1,\frac1{\rho L^2}\right\}I_{v,L}+J_L/\rho\right].
\]

For an arbitrary $T>0$, let $N=\max\{1,\lfloor\sqrt\rho\,T\rfloor\}$ and partition $[0,T]$ into $N$ equal intervals. Their common length $L=T/N$ lies in $(0,2\rho^{-1/2})$. Moreover, $\max\{1,(\rho L^2)^{-1}\}=\max\{1,(\rho T^2)^{-1}\}$. Theorem~\ref{thm:weak} gives finite entropy at each partition point, so the above estimate applies after a time translation on every interval. Averaging over the partition and taking the infimum over $W$ proves \eqref{eq:long-lsi}.
\end{proof}

\begin{remark}[Optimality of the short-time factor]\label{rem:lsi-short-time}
For the harmonic path of \cref{rem:endpoint} with $U(x)=\rho x^2/2$ and amplitude $r>0$, the intrinsic action vanishes and $\Ent_{\muT}(g)=\rho r^2/2$. Since $V_0=0$, the velocity Fisher information $\Iv(g_t)=V_t^2$ is of order $t^2$, and its average satisfies
\[
 I_{v,T}(g)=\frac{\rho r^2}{2}
       \left[1-\frac{\sin(2\sqrt\rho\,T)}{2\sqrt\rho\,T}\right]
       =\frac{\rho^2r^2}{3}T^2+O(T^4).
\]
Thus $\Ent_{\muT}(g)/I_{v,T}(g)=3/(2\rho T^2)+O(1)$ as $T\downarrow0$. The coefficient of $I_{v,T}$ in a boundary-free LSI must therefore grow at least as fast as $(\rho T^2)^{-1}$, even for a quadratic potential and paths of zero intrinsic action.
\end{remark}

For convex potentials with a spatial LSI, \cref{thm:long-lsi} extends the space-time LSI of~\cite{LiLu2026}*{Theorem 3.7 and Corollary 3.14}, with a time scale $T$ of order $\rho^{-1/2}$. It removes the Hessian-growth condition and the assumptions on the conditional mean velocity, higher velocity moments, and force moments imposed there. Admissible pairs in the sense of~\cite{LiLu2026}*{Definition 3.9} satisfy the technical assumptions of \cref{thm:weak}: the moment conditions (Q2) and (Q3) there and finite quadratic action give a $W_2$-continuous representative on $[0,T]$~\cite{ambrosio2005gradient}*{Theorem 8.3.1}, and the argument for \eqref{eq:truncatedchain}, started at a time of finite entropy, gives finite entropy at every time by lower semicontinuity. Hypercontractivity requires the additional interpolation argument of~\cite{LiLu2026}*{Sections 4--5}.

\section{Entropy decay for underdamped Langevin dynamics}\label{sec:decay}

In this section, we use the kinetic EVI \eqref{eq:intro-langevin-comparison} and entropy balance \eqref{eq:semigroup-entropy-balance} to prove entropy decay on every time interval under semiconvexity and a compatible spatial transportation inequality. This gives a cubic short-time bound and exponential decay for every friction coefficient, with the kinetic time scale obtained by choosing $\gamma=\sqrt\rho$. We first record a controlled estimate under positive curvature as a consequence of \cref{thm:weak}.

Throughout this section, we assume \eqref{eq:intro-potential} and the 
Talagrand's transportation inequality \eqref{eq:spatial-t2} with constant $\rho$. Under the positive curvature bound $\nabla^2U\succeq K\Id$, $K>0$, the latter follows with $\rho=K$ from \eqref{eq:classical-lsi} and the Otto--Villani theorem~\cite{OttoVillani2000}. By \eqref{eq:entropy-second-moment}, every density of finite entropy has a finite second moment. Proposition~\ref{prop:semigroup-paths} then shows the finite-action conditions and the entropy balance for its Langevin evolution, to which \cref{thm:kinetic-evi} applies.

\subsection{Entropy decay under positive curvature}

The following proposition is a variant of the modified entropy estimate of~\cite{Lu2026} and its controlled form in~\cite{LiLu2026}*{Proposition 3.6}. In the modified functional \eqref{eq:heavy-ball-functional} below, differentiating the current term $\beta\mathcal C_{\mux}$ and applying \eqref{eq:weak-comparison} with $W=gu-\gamma\nabla_vg$ yields spatial entropy dissipation together with the indefinite damping term $-\beta\gamma\mathcal C_{\mux}(g_t)$. The derivative of the additional distance term $\beta\gamma\mathcal D_{\mux}$ exactly cancels this damping term, as suggested by \eqref{eq:heavy-ball-intro} for the heavy-ball case. 

For $K>0$ and $\gamma>0$, set $\beta=c_0\min\{\gamma,K/\gamma\}$, where $c_0>0$ is a sufficiently small universal constant. For $\nu=g\mu$ of finite entropy, define
\begin{equation}\label{eq:heavy-ball-functional}
 \mathcal H_\beta(\nu)=\Ent_\mu(g)+\beta\mathcal C_{\mux}(g) +\beta\gamma\mathcal D_{\mux}(g).
\end{equation}

\begin{proposition}[Controlled modified entropy estimate]\label{prop:heavy-ball}
Assume $\nabla^2U\succeq K\Id$ for some $K>0$, and fix $\gamma>0$. Let $\nu_t=g_t\mu$ be narrowly continuous on $[0,T]$, with a density at every time, $\nu_0\in\PPtwo$, and $\Ent_\mu(g_0)<\infty$. Suppose that a measurable velocity-control field $u$ satisfies
\begin{equation}\label{eq:controlled-langevin}
 (\partial_t+\La)g=\gamma\Ls g+\nabla_v^*(gu),
 \qquad
 \int_0^T\left[\Iv(g_t)+\int g_t|u_t|^2\ud\mu\right]\ud t<\infty.
\end{equation}
Then $t\mapsto\mathcal H_\beta(\nu_t)$ is continuous and of bounded variation on $[0,T]$. It satisfies
\begin{equation}\label{eq:heavy-ball-equivalence}
 \frac12\Ent_\mu(g_t)
 \le\mathcal H_\beta(\nu_t)
 \le2\Ent_\mu(g_t),
\end{equation}
and, for a universal constant $C$, in $\mathcal D'(0,T)$,
\begin{equation}\label{eq:heavy-ball-controlled}
 \frac{\ud}{\ud t}\mathcal H_\beta(\nu_t)
       +\frac\beta2\mathcal H_\beta(\nu_t)
 \le\frac C\gamma  \int g_t|u_t|^2\ud\mu.
\end{equation}
In particular, for every $0\le s\le t\le T$,
\begin{align}\label{eq:heavy-ball-integrated}
 \mathcal H_\beta(\nu_t)\le{}&
 \e^{-\beta(t-s)/2}\mathcal H_\beta(\nu_s) + \frac C\gamma
       \int_s^t\e^{-\beta(t-r)/2}\!\int g_r|u_r|^2\ud\mu\ud r.
\end{align}
\end{proposition}

\begin{proof}
The path $g$ is driven by the current $W=gu-\gamma\nabla_vg$ that satisfies
\[
 \int\frac{|W_t|^2}{g_t}\ud\mu  \le2\int g_t|u_t|^2\ud\mu+2\gamma^2\Iv(g_t).
\]
Theorem~\ref{thm:weak} then gives the entropy chain rule, continuity, and bounded variation of $\mathcal C_{\mux}(g_t)$, and $\frac{\ud}{\ud t}\mathcal D_{\mux}(g_t)=\mathcal C_{\mux}(g_t)$. This proves the asserted regularity of $\mathcal H_\beta$. The transportation inequality with $\rho=K$ gives $K\mathcal D_{\mux}(g_t)\le\Ent_{\mux}(q_t)$, while \eqref{eq:current-t2} gives $|\mathcal C_{\mux}(g_t)|\le K^{-1/2}\Ent_\mu(g_t)$. Since $\beta\le c_0\sqrt K$ and $\beta\gamma\le c_0K$, these bounds imply \eqref{eq:heavy-ball-equivalence} for $c_0$ sufficiently small.

The entropy chain rule, \eqref{eq:weak-comparison}, \eqref{eq:score-mean}, and $\frac{\ud}{\ud t}\mathcal D_{\mux}(g_t)=\mathcal C_{\mux}(g_t)$ give, after Gaussian LSI,
\[
 \begin{aligned}
 \frac{\ud}{\ud t}\mathcal H_\beta(\nu_t) \le{}&-\beta\Ent_{\mux}(q_t)-\beta K\mathcal D_{\mux}(g_t)
       -(\gamma-\beta)\Iv(g_t)\\
 &+\int g_tu_t\cdot\nabla_v\log g_t\ud\mu +\beta\int g_tu_t\cdot\xi_{q_t}\ud\mu,
 \end{aligned}
\]
where the two terms associated with $u$ satisfy
\[
 \begin{aligned}
 \left|\int g_tu_t\cdot\nabla_v\log g_t\ud\mu\right|
 &\le\frac\gamma4\Iv(g_t)+\frac1\gamma\int g_t|u_t|^2\ud\mu,\\
 \beta\left|\int g_tu_t\cdot\xi_{q_t}\ud\mu\right|
 &\le\frac{\beta K}{2}\mathcal D_{\mux}(g_t) +\frac\beta K\int g_t|u_t|^2\ud\mu.
 \end{aligned}
\]
Using $\beta\le c_0\gamma$ and $\beta/K\le c_0/\gamma$, with $c_0$ sufficiently small, gives in $\mathcal D'(0,T)$,
\begin{align*}
 \frac{\ud}{\ud t}\mathcal H_\beta(\nu_t) \le{}&-\beta\Ent_{\mux}(q_t)-\frac{\beta K}{2}\mathcal D_{\mux}(g_t) -\frac\gamma2\Iv(g_t) +\frac C\gamma\int g_t|u_t|^2\ud\mu.
\end{align*}
Since $\beta\le\gamma$, Gaussian LSI and \eqref{eq:heavy-ball-equivalence} yield
\[
 \beta\Ent_{\mux}(q_t)
 +\frac\gamma2\Iv(g_t)
 \ge\beta\Ent_\mu(g_t)
 \ge\frac\beta2\mathcal H_\beta(\nu_t),
\]
which proves \eqref{eq:heavy-ball-controlled}. An integrating factor and continuity of $\mathcal H_\beta(\nu_t)$ give \eqref{eq:heavy-ball-integrated} for every $0\le s\le t\le T$.
\end{proof}

For $u=0$ and every initial density $g$ of finite entropy, \cref{prop:semigroup-paths} gives the path assumptions of \cref{prop:heavy-ball}, and then \eqref{eq:heavy-ball-equivalence} and \eqref{eq:heavy-ball-integrated} give, with universal constants $c,C>0$,
\begin{equation}\label{eq:positive-decay}
 \Ent_\mu(\Pt_tg) \le C\exp\!\left[-ct  \min\{\gamma,K/\gamma\}\right]\Ent_\mu(g), \qquad t\ge0.
\end{equation}
Compared with~\cite{Lu2026}*{Theorem 2.3}, the entropy decay \eqref{eq:positive-decay} assumes the stronger curvature bound $\nabla^2U\succeq K\Id$ with $K>0$: the earlier result assumes convexity and a spatial LSI with constant $\rho$. However, \eqref{eq:positive-decay} requires none of the polynomial growth and derivative bounds imposed in~\cite{Lu2026}*{Assumption 2.2}. Both estimates hold for every initial density of finite relative entropy and every $\gamma>0$. For uniformly convex potentials, taking $\rho=K$ in the earlier result gives the same rate order $\min\{\gamma,K/\gamma\}$ as \eqref{eq:positive-decay}.

\subsection{Entropy decay under semiconvexity}
Next, we prove the entropy contraction under semiconvexity and a compatible transportation inequality, on each interval of a prescribed length. We will quantify the contraction coefficient in terms of the interval length and the friction coefficient.
 
\begin{proposition}[Entropy decay on each time interval]\label{prop:entropy-drop}
Assume $\nabla^2U\succeq K\Id$ for some $K\le0$. Suppose that $\mux$ satisfies the transportation inequality \eqref{eq:spatial-t2} and that $\alpha=1+K/\rho>0$. There is a universal constant $c\in(0,1)$ such that, for every $\gamma>0$, every $\tau>0$, and every density $g$ of finite entropy,
\begin{equation}\label{eq:entropy-drop}
 \Ent_\mu(\Pt_\tau g)\le(1-\delta_\tau)\Ent_\mu(g),
\end{equation}
where
\begin{equation}\label{eq:delta}
 \delta_\tau=c\min\left\{1,\alpha\gamma\tau, \alpha^2\rho\gamma\tau^3,\frac{\alpha^2\rho\tau}{\gamma}\right\}.
\end{equation}
\end{proposition}

\begin{proof}
We write $g_t=\Pt_tg$, and below let $C$ denote a universal constant that may change from line to line. 
Substituting $W=-\gamma\nabla_vg$ into \eqref{eq:scalar-controlled} and using the score identity \eqref{eq:score-mean}, together with the entropy balance \eqref{eq:semigroup-entropy-balance}, we have the following estimates: 
\begin{equation} \label{eq:scalar-physical}
\begin{aligned}
\frac{\ud}{\ud t}\Ent_\mu(g_t)=-\gamma \Iv(g_t),\qquad |\mathcal C_{\mux}(g_t)|^2\le\frac2\rho\Ent_\mu(g_t)\Iv(g_t),\\
 \frac{\ud}{\ud t}\mathcal C_{\mux}(g_t)+\gamma\mathcal C_{\mux}(g_t)+\alpha\Ent_\mu(g_t)
 \le\frac{2+\alpha}{2}\Iv(g_t).        
\end{aligned}
\end{equation}
Here, as in \eqref{eq:scalar-controlled}, the differential inequality for $\mathcal C_{\mux}$ holds in distributions, and the other relations hold almost everywhere.

We test the differential inequality in \eqref{eq:scalar-physical} with the nonnegative cutoff $\chi_\tau$ in \eqref{eq:time-cutoff}. Since $\chi_\tau$ is compactly supported in $(0,\tau)$, integration by parts gives
\[ \alpha\int_0^\tau\chi_\tau\Ent_\mu(g_t)\ud t \le\frac{2+\alpha}{2}\int_0^\tau\chi_\tau\Iv(g_t)\ud t +\int_0^\tau\Bigl(\frac{\ud\chi_\tau}{\ud t}(t)-\gamma\chi_\tau\Bigr)\mathcal C_{\mux}(g_t)\ud t.
\]
The entropy $\Ent_\mu(g_t)$ is nonincreasing in $t$ and $\int_0^\tau\chi_\tau=1$, so the left-hand side is bounded below by $\alpha\Ent_\mu(g_\tau)$. On the right-hand side, $\alpha\le1$ and $\|\chi_\tau\|_\infty\le C/\tau$ from \eqref{eq:time-cutoff-norms} bound the first integral by $\frac C\tau\int_0^\tau\Iv(g_t)\ud t$. Hence,  we have 
\begin{equation} \label{eq:auxeqdecay1}
    \alpha\Ent_\mu(g_\tau)\le\frac C\tau\int_0^\tau \Iv(g_t)\ud t +\int_0^\tau\Bigl(\frac{\ud\chi_\tau}{\ud t}(t)-\gamma\chi_\tau\Bigr)\mathcal C_{\mux}(g_t)\ud t,
\end{equation}
where, by using \eqref{eq:time-cutoff-norms} and \eqref{eq:scalar-physical}, as well as $\Ent_\mu(g_t)\le\Ent_\mu(g)$, 
\begin{equation} \label{eq:auxeqdecay2}
      \begin{aligned}
\left|\int_0^\tau\Bigl(\frac{\ud\chi_\tau}{\ud t}(t)-\gamma\chi_\tau\Bigr)\mathcal C_{\mux}(g_t)\ud t\right| & \le\left(\int_0^\tau\Bigl(\frac{\ud\chi_\tau}{\ud t}(t)-\gamma\chi_\tau\Bigr)^2\ud t\right)^{1/2}
 \left(\int_0^\tau|\mathcal C_{\mux}(g_t)|^2\ud t\right)^{1/2}\\
 & \le C\left(\frac1{\tau^3}+\frac{\gamma^2}\tau\right)^{1/2}
 \Bigl(\frac{\Ent_\mu(g)}{\rho}\int_0^\tau \Iv(g_t)\ud t\Bigr)^{1/2}.
 \end{aligned}
\end{equation}

We write the decay of the entropy over $[0,\tau]$ as
\[
 \Delta_\tau=\Ent_\mu(g)-\Ent_\mu(g_\tau) =\gamma\int_0^\tau\Iv(g_t)\ud t.
\]
The above estimates \eqref{eq:auxeqdecay1}-\eqref{eq:auxeqdecay2} and Young's inequality give
\[
 \begin{aligned}
 \alpha\Ent_\mu(g)
 &\le\left(\alpha+\frac C{\gamma\tau}\right)\Delta_\tau
       +C\left[\frac{\Ent_\mu(g)\Delta_\tau}{\rho\gamma}
          \left(\frac1{\tau^3}+\frac{\gamma^2}\tau\right)\right]^{1/2}\\
 &\le\frac\alpha2\Ent_\mu(g)
       +\left[\alpha+\frac C{\gamma\tau}
          +\frac C{\alpha\rho\gamma}
             \left(\frac1{\tau^3}+\frac{\gamma^2}\tau\right)\right]\Delta_\tau.
 \end{aligned}
\]
It follows that 
\[
 \Ent_\mu(g)\le C\left[ 1+\frac1{\alpha\gamma\tau}  +\frac1{\alpha^2\rho\gamma\tau^3} +\frac{\gamma}{\alpha^2\rho\tau}
 \right]\Delta_\tau,
\]
which proves \eqref{eq:entropy-drop} by taking $c>0$ sufficiently small in \eqref{eq:delta}.
\end{proof}

\begin{corollary}[Cubic short-time bound and exponential decay]\label{cor:decay}
Assume $\nabla^2U\succeq K\Id$ with $K\le0$, the transportation inequality \eqref{eq:spatial-t2}, and $\alpha=1+K/\rho>0$. There is a universal constant $c>0$ such that, for every $\gamma>0$ and every density $g$ of finite entropy,
\begin{equation}\label{eq:decay}
 \Ent_\mu(\Pt_tg)
 \le\exp\!\left[-c\,\alpha^2
          \min\{\rho\gamma t^3,\gamma t,\rho t/\gamma\}\right]\Ent_\mu(g),
 \qquad t\ge0.
\end{equation}
In particular, the choice $\gamma=\sqrt\rho$ gives
\begin{equation}\label{eq:decay-kinetic}
 \Ent_\mu(\Pt_tg)
 \le\exp\!\left[-c\,\alpha^2
          \min\{(\sqrt\rho\,t)^3,\sqrt\rho\,t\}\right]\Ent_\mu(g).
\end{equation}
For convex potentials one may take $K=0$ and $\alpha=1$.
\end{corollary}

\begin{proof}
Let $\lambda=\min\{\gamma,\rho/\gamma\}$. Then $\lambda\le\gamma$, $\gamma\lambda\le\rho$, and $\lambda^3\le\rho\gamma$. Throughout the proof, $c>0$ is a universal constant that may decrease from line to line. For $0<t\le\lambda^{-1}$, \eqref{eq:delta} and $0<\alpha\le1$ give
\[
 \delta_t\ge c\alpha^2\min\{\rho\gamma t^3,\lambda t\}.
\]
Proposition~\ref{prop:entropy-drop} and $1-r\le\e^{-r}$ therefore prove \eqref{eq:decay} for $0<t\le\lambda^{-1}$.

For $t\ge\lambda^{-1}$, let $N=\lfloor\lambda t\rfloor$. At the fixed interval length $\lambda^{-1}$, the contraction \eqref{eq:delta} gives $\delta_{\lambda^{-1}}\ge c\alpha^2$. Iterating \eqref{eq:entropy-drop} over $N$ such intervals and using entropy monotonicity on the remaining interval yields
\[
 \Ent_\mu(\Pt_tg)
 \le\Ent_\mu(\Pt_{N/\lambda}g)
 \le(1-\delta_{\lambda^{-1}})^N\Ent_\mu(g)
 \le\exp\!\left[-c\alpha^2\lambda t\right]\Ent_\mu(g),
\]
where we used $N\ge\lambda t/2$ in the last step. For $t\ge\lambda^{-1}$, the inequality $\lambda^3\le\rho\gamma$ implies $\rho\gamma t^3\ge\lambda t$. This proves \eqref{eq:decay} in the remaining time range. The case $t=0$ is immediate, and \eqref{eq:decay-kinetic} follows by setting $\gamma=\sqrt\rho$. We hence complete the proof. 
\end{proof}

\begin{remark}[Optimality of the cubic short-time order]\label{rem:cubic-sharp}
Similarly to \cref{rem:lsi-short-time} and \cref{rem:endpoint}, we take the potential $U(x)=\rho x^2/2$ in one dimension and the initial distribution $g\mu=\mathcal N(r,\rho^{-1})\otimes\mathcal N(0,1)$ with $r\ne0$. The Langevin evolution is
\[
 (\Pt_tg)\mu=\mathcal N(X_t,\rho^{-1})\otimes\mathcal N(V_t,1),
 \qquad
 \frac{\ud X_t}{\ud t}=V_t,\quad \frac{\ud V_t}{\ud t}=-\rho X_t-\gamma V_t,
\]
with $(X_0,V_0)=(r,0)$. The covariance of $(\Pt_tg)\mu$ is constant in time and hence remains equal to that of $\mu$. In addition, we compute
\[
 \Ent_\mu(\Pt_tg)=\tfrac12(\rho X_t^2+V_t^2),\qquad
 \frac{\ud}{\ud t}\Ent_\mu(\Pt_tg)=-\gamma V_t^2.
\]
Since $V_t=-\rho r t+O(t^2)$ and $\Ent_\mu(g)=\frac{\rho r^2}{2}>0$, we integrate $\frac{\ud}{\ud t}\Ent_\mu(\Pt_tg)$ and find
\begin{equation}\label{eq:cubic-sharp}
 \frac{\Ent_\mu(\Pt_tg)}{\Ent_\mu(g)}
 =1-\frac{2\rho\gamma}{3}t^3+O(t^4),\qquad t\downarrow0.
\end{equation}
Thus the cubic power at short time in \eqref{eq:decay} is optimal for entropy decay bounds with prefactor one over all finite-entropy initial densities.
\end{remark}

Under the positive curvature assumption $\nabla^2U\succeq K\Id$ with $K>0$, we may apply \cref{cor:decay} using the weaker Hessian bound $\nabla^2U\succeq0$, together with $\rho=K$ and $\alpha=1$. The resulting estimate recovers the long-time decay rate of \eqref{eq:positive-decay} up to universal constants and also provides the cubic short-time bound. Its prefactor is $1$ because \eqref{eq:entropy-drop} compares $\Ent_\mu(\Pt_\tau g)$ directly with $\Ent_\mu(g)$.

\appendix
\section{Entropy chain rule and approximation}\label{app:recovery}

This appendix proves the entropy chain rule, the superposition and approximation results, and the transport stability used in \cref{thm:weak}, together with \cref{prop:semigroup-paths}. We follow the density-current formulation of~\cite{LiLu2026}*{Appendix B}.

Throughout this appendix, we assume \eqref{eq:intro-potential}, and write $z=(x,v)$, $\varrho_t=g_t\varrho_\infty$, and $F_t=W_t\varrho_\infty=\varrho_tA_t$, as in \eqref{eq:ordinary-control}. Derivatives are understood in distributions, with the zero-density conventions of \eqref{eq:intro-path}.

\subsection{The entropy chain rule}\label{app:entropy-chain}

\begin{proposition}[Entropy chain rule]\label{prop:entropy-chain}
Under the assumptions of \cref{thm:weak}, the entropy $\Ent_\mu(g_t)$ is finite and absolutely continuous on $[0,T]$, and \eqref{eq:entropy-chain} holds for every $0\le s\le t\le T$.
\end{proposition}

The following lemma supplies the strong convergence needed to pass to the limit in the product of the velocity current and the velocity score.

\begin{lemma}[Convergence of weighted fields]\label{lem:weighted-fields}
Let $r_k,r$ be nonnegative functions and $B_k,B$ vector fields on a common open subset of Euclidean space. Suppose $r_k\to r$ in $L^1$, $B_k\to B$ locally in distributions, $B_k=0$ almost everywhere on $\{r_k=0\}$, $B=0$ almost everywhere on $\{r=0\}$, and
\[
 \int|B_k|^2/r_k\longrightarrow\int|B|^2/r<\infty.
\]
Then $B_k/\sqrt{r_k}\to B/\sqrt r$ strongly in $L^2$.
\end{lemma}

\begin{proof}
Noting that the quotients $B_k/\sqrt{r_k}$ are bounded in $L^2$, one can choose a subsequence with weak limit $b$. Since $r_k\to r$ in $L^1$, we have $\sqrt{r_k}\to\sqrt r$ in $L^2$. Testing the products $\sqrt{r_k}(B_k/\sqrt{r_k})$ against smooth compactly supported functions therefore gives $\sqrt r\,b=B$, so $b=B/\sqrt r$ on $\{r>0\}$.

Weak lower semicontinuity and the assumed convergence of the squared norms give
\[
 \int\frac{|B|^2}{r}
 \le\|b\|_{L^2}^2
 \le\liminf_k\int\frac{|B_k|^2}{r_k}
 =\int\frac{|B|^2}{r},
\]
which implies that $b=0$ on $\{r=0\}$, and the norms converge to $\|b\|_{L^2}$. Hence, the weak convergence can be upgraded to the strong one.  Every subsequence has a further subsequence with this same strong limit, which proves convergence of the full sequence.
\end{proof}

In the rest of this subsection, let $(g,W)$ satisfy the assumptions of \cref{thm:weak}. 
For the cutoff argument, we use the modified Hamiltonian energy
\begin{equation}\label{eq:energy}
 E(x,v)=U(x)+c|x|^2+C_0+\tfrac12|v|^2
 \ge 1+|x|^2+\tfrac12|v|^2,
\end{equation}
where $c\ge0$ and $C_0$ are chosen sufficiently large such that the inequality holds. Such a choice is possible because semiconvexity gives a quadratic lower bound on $U$. Thus $E$ is smooth and proper. The force terms cancel in $\La E$, and
\begin{equation}\label{eq:energyderivatives}
 \La E=2c\,x\cdot v,\qquad
 \nabla_vE=v,\qquad |\La E|\le CE.
\end{equation}
For each fixed $\chi\in C_c^\infty([0,\infty))$, \eqref{eq:energyderivatives} implies that the cutoff derivatives $\La[\chi(E/R)]$ are bounded uniformly in $(x,v)\in\RR^{2d}$ and $R\ge1$. 

Let $S_t$ denote the velocity score weighted by the Lebesgue density:
\begin{equation}\label{eq:ordinary-score}
 S_t=\varrho_t\nabla_v\log g_t
     =\nabla_v\varrho_t+v\varrho_t.
\end{equation}
Then we have 
\[
 \int\frac{|F_t|^2}{\varrho_t}\ud z =\int\frac{|W_t|^2}{g_t}\ud\mu,
 \qquad
 \int\frac{|S_t|^2}{\varrho_t}\ud z=\Iv(g_t).
\]

We smooth the density-current pair in phase space using Galilean translations, as in~\cite{BrigatiMaasQuattrocchi2025}*{Section 4.3, (90)--(94)}, without averaging in time. Choose a centered, nonnegative $\theta\in C_c^\infty(\RR^{2d})$ of integral one and define the following smoothed functions: 
\begin{equation}\label{eq:shear-recovery}
 \begin{aligned}
 \theta_{\varepsilon,t}(x,v)
 &=\varepsilon^{-2d}
   \theta\bigl((x-tv)/\varepsilon,v/\varepsilon\bigr),
 \qquad
 \varrho_t^\varepsilon=\theta_{\varepsilon,t}*\varrho_t,\\
 F_t^\varepsilon
 &=\theta_{\varepsilon,t}*F_t
   +\nabla U\,\varrho_t^\varepsilon
   -\theta_{\varepsilon,t}*(\nabla U\,\varrho_t),\\
 S_t^\varepsilon
 &=\nabla_v\varrho_t^\varepsilon+v\varrho_t^\varepsilon.
 \end{aligned}
\end{equation}
All convolutions are taken in $(x,v)$ and are locally well-defined because the kernel is compactly supported. For a smooth function $f=f(t,x,v)$, we compute
\[
 \begin{aligned}
 &(\partial_t+v\cdot\nabla_x)(\theta_{\varepsilon,t}*f)
 -\theta_{\varepsilon,t}*((\partial_t+v\cdot\nabla_x)f)=
\bigl[(\partial_t+v\cdot\nabla_x)\theta_{\varepsilon,t}\bigr]*f=0. 
 \end{aligned}
\]
This identity extends to distributions. Thus, unlike the standard Gaussian convolution used in~\cite{LiLu2026}*{Appendix B}, this smoothing commutes with $\partial_t+v\cdot\nabla_x$. The correction $\nabla U\,\varrho_t^\varepsilon -\theta_{\varepsilon,t}*(\nabla U\,\varrho_t)$ in $F_t^\varepsilon$ cancels the commutator with the force operator $-\nabla U\cdot\nabla_v$.
Hence $(\varrho^\varepsilon,F^\varepsilon)$ solves
\eqref{eq:ordinary-control}.

\begin{lemma}[Local strong convergence]\label{lem:local-fields}
For every compact $Q\subset\RR^{2d}$,
\begin{equation}\label{eq:local-weighted-fields}
 \frac{F^\varepsilon}{\sqrt{\varrho^\varepsilon}}
 \longrightarrow\frac F{\sqrt\varrho},\qquad
 \frac{S^\varepsilon}{\sqrt{\varrho^\varepsilon}}
 \longrightarrow\frac S{\sqrt\varrho}
 \quad\text{strongly in }L^2((0,T)\times Q).
\end{equation}
Also $\varrho^\varepsilon\to\varrho$ in local space-time $L^1$ and in local $L^1$ at every fixed time, including $0$ and $T$.
\end{lemma}
\begin{proof}
We fix a compact set $Q\subset\RR^{2d}$ and choose a nonnegative $\psi\in C_c^\infty(\RR^{2d})$ with $\psi=1$ on $Q$. Choose a fixed closed ball $B$ centered at the origin such that $\operatorname{supp}\theta_{\varepsilon,t}\subset B$ for all $0<\varepsilon\le1$ and $t\in[0,T]$. Evaluating the convolution on $\operatorname{supp}\psi$ then involves only points in the compact set $\operatorname{supp}\psi-B$. The force $\nabla U$ is Lipschitz on the spatial projection of this set.  
For almost every $t$, let $Y$ be a random vector with density $\varrho_t$, and let $\eta$ be an independent random vector with density $\theta_{\varepsilon,t}$. Letting $Z = Y + \eta$, we have 
\begin{align*}
 \frac{F_t^\varepsilon}{\varrho_t^\varepsilon}(Z)
 &=\mathbb E[A_t(Y)+\nabla U(Y_x+\eta_x)-\nabla U(Y_x)\mid Z],\\
 \frac{S_t^\varepsilon}{\varrho_t^\varepsilon}(Z)
 &=\mathbb E[S_t(Y)/\varrho_t(Y)+\eta_v\mid Z].
\end{align*}
Both $\eta_v$ and the increment of $\nabla U$: $\nabla U(Y_x+\eta_x)-\nabla U(Y_x)$, are bounded by $C\varepsilon$ when $Z\in\operatorname{supp}\psi$. Using Jensen's inequality for the conditional expectation, weighted by $\psi(Z)$, and the triangle inequality, we find
\begin{equation}\label{eq:jensen-smoothing}
 \left(\int\psi\frac{|F_t^\varepsilon|^2}{\varrho_t^\varepsilon}\ud z\right)^{1/2}
 \le\left(\int\!\!\int\psi(z+y)\theta_{\varepsilon,t}(y)\ud y\,
        \frac{|F_t(z)|^2}{\varrho_t(z)}\ud z\right)^{1/2}+C\varepsilon,
\end{equation}
and the same bound with $S$ in place of $F$.  The weights $\int\psi(z+y)\theta_{\varepsilon,t}(y)\ud y$ converge to $\psi(z)$ uniformly in $z\in\RR^{2d}$ and $t\in[0,T]$ as $\varepsilon\downarrow0$, and they are supported in a fixed compact set. Integrating \eqref{eq:jensen-smoothing} in time therefore gives
\[
 \limsup_{\varepsilon\downarrow0}
 \int_0^T\!\int\psi\frac{|F^\varepsilon|^2}{\varrho^\varepsilon}\ud z\ud t
 \le\int_0^T\!\int\psi\frac{|F|^2}{\varrho}\ud z\ud t,
\]
and the same bound with $S$ in place of $F$.

The approximate-identity property gives local $L^1$ convergence of $\varrho^\varepsilon$ on space-time and at every fixed time. Also $F^\varepsilon\to F$ and $S^\varepsilon\to S$ locally in distributions. Lower semicontinuity therefore gives the reverse inequalities for the weighted quadratic integrals. Applying \cref{lem:weighted-fields} with $r_\varepsilon=\psi\varrho^\varepsilon$ and $B_\varepsilon=\psi F^\varepsilon$, and then with $B_\varepsilon=\psi S^\varepsilon$, we complete the proof of \eqref{eq:local-weighted-fields} by $\psi=1$ on $Q$. 
\end{proof}

We truncate the second derivative of $r\log r-r+1$ to derive entropy dissipation identities without assuming finite entropy at positive times (which is itself part of the conclusion of \cref{prop:entropy-chain}). Choose smooth functions $a_n:(0,\infty)\to[0,1]$, compactly supported in $(0,\infty)$, such that $a_n\uparrow1$, and define $\beta_n$ by
\begin{equation}\label{eq:beta}
 \frac{\ud^2\beta_n}{\ud r^2}(r)=\frac{a_n(r)}r,\qquad
 \beta_n(1)=\frac{\ud\beta_n}{\ud r}(1)=0.
\end{equation}
Extend $a_n$ by zero at $r=0$ and extend $\beta_n$ smoothly to $r=0$. Then $\beta_n\ge0$,
\[
 \beta_n(r)\uparrow r\log r-r+1,\qquad
 \left|\frac{\ud\beta_n}{\ud r}(r)\right|\le C_n,\qquad
 \beta_n(r)\le C_n(r+1),\qquad 0\le r\frac{\ud^2\beta_n}{\ud r^2}(r)\le1.
\]

\begin{lemma}[Local renormalization, including endpoint traces]\label{lem:local-chain}
For $\psi\in C_c^\infty(\RR^{2d})$ and $0\le s\le t\le T$,
\begin{align}\label{eq:localchain}
 \int\psi\beta_n(g_t)\ud\mu-\int\psi\beta_n(g_s)\ud\mu
 =\int_s^t\!\int\bigl[
 \beta_n(g)\La\psi+\frac{\ud\beta_n}{\ud r}(g)W\cdot\nabla_v\psi
 +\psi\frac{\ud^2\beta_n}{\ud r^2}(g)W\cdot\nabla_vg\bigr]\ud\mu\ud\tau.
\end{align}
\end{lemma}
\begin{proof}
We set $g^\varepsilon=\varrho^\varepsilon/\varrho_\infty$ and $W^\varepsilon=F^\varepsilon/\varrho_\infty$, which satisfies  $(\partial_t+\La)g^\varepsilon=\nabla_v^*W^\varepsilon$.
The density $g^\varepsilon$ is smooth in phase space, and its time derivative is integrable on $(0,T)\times Q$ for every compact $Q\subset\RR^{2d}$.
Narrow continuity of $\nu_t$ ensures that $g_t^\varepsilon$ is continuous in time, locally uniformly in phase space, so its time traces agree with the prescribed smoothed densities at every $t\in[0,T]$.
Since $\frac{\ud\beta_n}{\ud r}$ and $\frac{\ud^2\beta_n}{\ud r^2}$ are bounded, the chain rule for absolutely continuous functions yields \eqref{eq:localchain} for $(g^\varepsilon,W^\varepsilon)$ between any two times $0\le s\le t\le T$. We next pass to the limit as $\varepsilon\downarrow0$. 

Indeed, \cref{lem:local-fields} gives $\varrho^\varepsilon\to\varrho$ in local $L^1$. The factorization
$F^\varepsilon=\sqrt{\varrho^\varepsilon}\,(F^\varepsilon/\sqrt{\varrho^\varepsilon})$
and Cauchy--Schwarz then imply $F^\varepsilon\to F$ in local $L^1$.
Since $\beta_n$ is Lipschitz and $\frac{\ud\beta_n}{\ud r}$ is bounded and continuous, the terms containing $\La\psi$ and $\nabla_v\psi$ converge in $L^1$. In Lebesgue variables, the term containing $\frac{\ud^2\beta_n}{\ud r^2}$ has the following integrand:
\begin{equation} \label{eq:int}
     \psi\,\frac{\ud^2\beta_n}{\ud r^2}(g^\varepsilon) W^\varepsilon\cdot\nabla_vg^\varepsilon\,\varrho_\infty = \psi\,a_n(g^\varepsilon) \frac{F^\varepsilon}{\sqrt{\varrho^\varepsilon}} \cdot\frac{S^\varepsilon}{\sqrt{\varrho^\varepsilon}}.
\end{equation}
By \cref{lem:local-fields}, the two weighted fields converge strongly in local $L^2$, so their product converges in local $L^1$.
Together with $a_n(g^\varepsilon)\to a_n(g)$ in measure and $0\le a_n\le1$, this yields $L^1$ convergence of \eqref{eq:int} as $\varepsilon \downarrow 0$.
Finally, local $L^1$ convergence of the densities at times $s$ and $t$, combined with the Lipschitz continuity of $\beta_n$, yields convergence of both trace terms.
This proves \eqref{eq:localchain}. 
\end{proof}

\begin{proof}[Proof of \cref{prop:entropy-chain}]
We first remove the energy cutoff at fixed $n$ and then let $n\to\infty$. Choose a nonincreasing $\chi\in C_c^\infty([0,\infty))$ taking values in $[0,1]$, equal to one on $[0,1]$, and vanishing on $[2,\infty)$. Apply \cref{lem:local-chain} with $\psi_R=\chi(E/R)$. The identities \eqref{eq:energyderivatives} give
\begin{equation}\label{eq:cutoffbounds}
 |\La\psi_R|\le C\mathbf1_{\{E\ge R\}},\qquad
 |\nabla_v\psi_R|\le CR^{-1/2}\mathbf1_{\{E\ge R\}}.
\end{equation}
Since $\beta_n(g)\le C_n(g+1)$, the term containing $\La\psi_R$ is bounded in absolute value by
\[
 C_n\int_s^t\!\int_{\{E\ge R\}}(g+1)\ud\mu\ud\tau\longrightarrow0.
\]
Since $|\frac{\ud\beta_n}{\ud r}|\le C_n$, the term containing $W\cdot\nabla_v\psi_R$ is bounded in absolute value by
\[
 C_nR^{-1/2}(t-s)^{1/2}
 \left(\int_s^t\!\int\frac{|W|^2}{g}\ud\mu\ud\tau\right)^{1/2}
 \longrightarrow0.
\]
In the term containing $\psi_R\frac{\ud^2\beta_n}{\ud r^2}(g)W\cdot\nabla_vg$, the identity $g\frac{\ud^2\beta_n}{\ud r^2}(g)=a_n(g)\le1$ bounds the integrand by
\[
 g|A|\,|\nabla_v\log g|,
\]
which is integrable by Cauchy--Schwarz and \eqref{eq:intro-path}. The traces converge by dominated convergence, since $\beta_n(g)\le C_n(g+1)$. Letting $R\to\infty$ therefore yields
\begin{equation}\label{eq:truncatedchain}
 \int\beta_n(g_t)\ud\mu-\int\beta_n(g_s)\ud\mu
 =\int_s^t\!\int a_n(g)\,gA\cdot\nabla_v\log g\ud\mu\ud\tau.
\end{equation}
To remove the entropy truncation, take $s=0$ and use $\beta_n(g_0)\le g_0\log g_0-g_0+1$. Cauchy--Schwarz gives, uniformly in $n$ and $t$,
\begin{equation}\label{eq:entropybound}
 \int\beta_n(g_t)\ud\mu
 \le\Ent_\mu(g_0)+
 \left(\int_0^T\!\int\frac{|W|^2}{g}\ud\mu\ud\tau\right)^{1/2}
 \left(\int_0^T \Iv(g_\tau)\ud\tau\right)^{1/2}.
\end{equation}
Monotone convergence proves that the entropy is finite at every time. We may now let $n\to\infty$ in \eqref{eq:truncatedchain}: monotone convergence applies to both traces and dominated convergence to the time integral. This gives \eqref{eq:entropy-chain} for every $0\le s\le t\le T$. The integrability of $gA\cdot\nabla_v\log g$ proves absolute continuity.
\end{proof}

\subsection{Superposition of controlled trajectories}\label{app:superposition}

We represent controlled paths by trajectories and bound their maximal second moment. This also proves continuity of the path in quadratic Wasserstein distance.

\begin{lemma}[Superposition and a pathwise moment bound]\label{lem:superposition}
Assume \eqref{eq:intro-potential}. Let $\nu_t=g_t\mu$ be narrowly continuous on $[0,T]$, have a density at every time, and solve \eqref{eq:ordinary-control} with $F=\varrho A$. Suppose that $\nu_0\in\PPtwo$, $\Ent_\mu(g_0)<\infty$, and $\int_0^T\int|A|^2\ud\nu_t\ud t<\infty$. Then $\nu\in C([0,T];\PPtwo)$, and there is a probability measure $\mathbf P$ on absolutely continuous paths $Z=(X,V)$ such that
\begin{equation}\label{eq:characteristics}
 \dot X_t=V_t,\qquad
 \dot V_t=-\nabla U(X_t)+A_t(Z_t),\qquad (Z_t)_\#\mathbf P=\nu_t.
\end{equation}
Moreover,
\begin{equation}\label{eq:supmoment}
 \mathbb E_{\mathbf P}\sup_{0\le t\le T}|Z_t|^2
 \le C_T\left(\int E\ud\nu_0+
 \int_0^T\!\int|A|^2\ud\nu_t\ud t\right)<\infty.
\end{equation}
\end{lemma}

We first show that the initial modified Hamiltonian energy is integrable.

\begin{lemma}[Energy integrability]\label{lem:energy-integrability}
If $\nu=g\mu\in\PPtwo$ and $\Ent_\mu(g)<\infty$, then $\int E\ud\nu<\infty$.
\end{lemma}
\begin{proof}
By the semiconvexity of the potential, the negative part of $U$ is bounded by $C(1+|x|^2)$. Thus $f(x,v)=U(x)-|x|^2$ has a $\nu$-integrable negative part, and
\[
 \int \e^f\ud\mu=Z_x^{-1}\int \e^{-|x|^2}\ud x<\infty.
\]
Applying the entropy variational inequality to upper and lower truncations of $f$ gives
\[
 \int U\ud\nu\le\Ent_\mu(g)+\int|x|^2\ud\nu
 +\log\left(Z_x^{-1}\int \e^{-|x|^2}\ud x\right).
\]
Thus, $U^+$ is $\nu$-integrable. The definition \eqref{eq:energy} and the finite second moment of $\nu$ now give $\int E\ud\nu<\infty$.
\end{proof}

\begin{proof}[Proof of \cref{lem:superposition}]
We apply the superposition theorem of~\cite{BogachevRocknerShaposhnikov2019}*{Theorem 1.1}, whose global integrability condition involves only the radial component of the drift. We embed phase space so that the squared radius equals $4E$. The radial component is then expressed in terms of the derivative of $E$ along the drift, in which the terms involving $\nabla U$ cancel. 

We choose a Borel representative of $A$ and set $b(t,z)=(v,-\nabla U(x)+A_t(z))$. Define the smooth proper embedding
\[
 \iota(z)=\left(z,\sqrt{4E(z)-|z|^2}\right)\in\RR^{2d+1}.
\]
The lower bound in \eqref{eq:energy} ensures that the expression under the square root is strictly positive. Moreover, $|\iota(z)|^2=4E(z)$. We further set $\widehat\nu_t=\iota_\#\nu_t$ and define the lifted drift by
\[
 \widehat b(t,\iota(z))=D\iota(z)b(t,z),
\]
extending $\widehat b$ by zero outside the closed graph $\iota(\RR^{2d})$.
The path $t\mapsto\widehat\nu_t$ is narrowly continuous, since $\iota$ is continuous. 
It also satisfies the continuity equation with drift $\widehat b$.
Indeed, for every $\phi\in C_c^\infty(\RR^{2d+1})$, smoothness and properness of $\iota$ give $\phi\circ\iota\in C_c^\infty(\RR^{2d})$.
Testing the original equation with $\phi\circ\iota$ and applying the chain rule yields the desired equation.

By the finite action and Cauchy--Schwarz inequality, we have $A\in L^1(\ud t\,\nu_t)$. Since $\nabla U$ is locally bounded and $\iota$ is smooth and proper, the lifted drift is locally integrable with respect to $\ud t\,\widehat\nu_t$. Moreover,
\[
 \widehat b(t,\iota(z))\cdot\iota(z)
 =\tfrac12 b(t,z)\cdot\nabla(4E(z))
 =4c\,x\cdot v+2A_t(z)\cdot v.
\]
The lower bound in \eqref{eq:energy} therefore gives
\begin{equation}\label{eq:radial}
 \frac{|\widehat b(t,\iota(z))\cdot\iota(z)|}
 {(1+|\iota(z)|)^2}\le C(1+|A_t(z)|).
\end{equation}
The right-hand side is integrable with respect to $\ud t\,\nu_t$. Applying the superposition theorem with zero diffusion matrix thus yields a probability measure on $C([0,T];\RR^{2d+1})$ with time marginals $\widehat\nu_t$. 

We next show that these paths solve the ordinary differential equation driven by $\widehat b$.
Let $\widehat Z$ denote the coordinate process and set
$\tau_R=\inf\{t\ge0:|\widehat Z_t|\ge R\}\wedge T$.
Local integrability of $\widehat b$, together with the martingale identities for smooth cutoffs of the coordinate functions and their products, yields continuous martingales, with vanishing quadratic covariations from the zero diffusion matrix,
\[
 \begin{aligned}
 M_t^{R,i}
 &=\widehat Z_{t\wedge\tau_R}^{i}-\widehat Z_0^{i}
   -\int_0^{t\wedge\tau_R}\widehat b^{i}(s,\widehat Z_s)\ud s,\\
 \langle M^{R,i},M^{R,j}\rangle_t&=0,
 \qquad 1\le i,j\le2d+1.
 \end{aligned}
\]
Since $M_0^{R,i}=0$, we have $M^{R,i}\equiv0$ almost surely.
Thus the paths are absolutely continuous up to $\tau_R$ and satisfy
$\dot{\widehat Z}_t=\widehat b(t,\widehat Z_t)$ almost everywhere before $\tau_R$.
Every continuous path on $[0,T]$ is bounded, so letting $R\to\infty$ removes the stopping. The prescribed marginals $\widehat\nu_t$ imply that, almost surely, the paths lie on $\iota(\RR^{2d})$ at every rational time.
Since this graph is closed and the paths are continuous, they lie on it at every time.
Projecting onto the first $2d$ coordinates gives \eqref{eq:characteristics}.  

Finally, the marginal identity and finite action imply $\int_0^T|A_t(Z_t)|^2\ud t<\infty$ almost surely. Along each such trajectory, we have 
\[
 \frac{\ud}{\ud t}E(Z_t)
 =2cX_t\cdot V_t+A_t(Z_t)\cdot V_t
 \le CE(Z_t)+\tfrac12|A_t(Z_t)|^2.
\]
Gronwall's inequality bounds $\sup_{t\le T}E(Z_t)$ by $C_T\bigl(E(Z_0)+\int_0^T|A_t(Z_t)|^2\ud t\bigr)$. Taking expectation, using \cref{lem:energy-integrability}, and recalling $|Z_t|^2\le2E(Z_t)$ proves \eqref{eq:supmoment}.

In particular, $\nu_t\in\PPtwo$ for every $t\in[0,T]$. Coupling $\nu_s$ and $\nu_t$ by the same trajectory gives
\[
 W_2^2(\nu_s,\nu_t)\le\mathbb E_{\mathbf P}|Z_s-Z_t|^2
 \longrightarrow0\qquad\text{as }s\to t.
\]
Indeed, the trajectories are continuous and $|Z_s-Z_t|^2\le4\sup_{r\le T}|Z_r|^2$. Hence, \eqref{eq:supmoment} permits dominated convergence, including at the endpoints, which proves $\nu\in C([0,T];\PPtwo)$. 
\end{proof}

\subsection{Langevin paths and entropy balance}\label{app:langevin-paths}

\begin{proof}[Proof of \cref{prop:semigroup-paths}]
We first establish nonexplosion and invariance, then the entropy balance, and finally $W_2$ continuity.
Let $E$ be the modified Hamiltonian energy in \eqref{eq:energy}. For the diffusion generator $L=\La+\gamma\Ls$ acting on observables,
\[
 LE=2c\,x\cdot v-\gamma|v|^2+\gamma d,\qquad |LE|\le CE.
\]
The coefficients are smooth, so the Langevin equation has a unique solution up to explosion. For $R>E(z)$, set $\tau_R=\inf\{t\ge0:E(Z_t)\ge R\}$. It\^o's formula gives
\[
 \begin{aligned}
 E(Z_{t\wedge\tau_R})
 &=E(z)+\int_0^{t\wedge\tau_R}LE(Z_s)\ud s+M_t^R,\\
 M_t^R&=\sqrt{2\gamma}\int_0^{t\wedge\tau_R}V_s\cdot\ud B_s.
 \end{aligned}
\]
The stopped martingale $M^R$ has mean zero, so taking expectations and using Gronwall's inequality yields
\[
 \mathbb E_z E(Z_{t\wedge\tau_R})\le\e^{Ct}E(z),
 \qquad
 \mathbb P_z(\tau_R\le T)\le\frac{\e^{CT}E(z)}R\longrightarrow0.
\]
Properness of $E$ therefore gives nonexplosion.

To prove invariance, first note that the constant path $\mu$ solves the stationary Fokker--Planck equation against compactly supported smooth tests. Use the embedding $\iota$ from \cref{lem:superposition} to lift this path. Set $w=\sqrt{4E-|z|^2}$. By \eqref{eq:energy},
\[
 w^2\ge1+|z|^2,\qquad \nabla_vw=v/w.
\]
Thus $D_v\iota$ is bounded, so the lifted diffusion matrix below is bounded. The lifted coefficients are
\[
 \begin{aligned}
 \widehat a&=\gamma D_v\iota(D_v\iota)^\mathsf T,\\
 \widehat b&=D\iota\,(v,-\nabla U-\gamma v)+\gamma\Delta_v\iota.
 \end{aligned}
\]
Applying $L$ to $|\iota|^2=4E$ yields
\[
 2\widehat b\cdot\iota+2\tr\widehat a
 =L|\iota|^2=4LE.
\]
Since $|LE|\le CE$ and $|\iota|^2=4E$, the radial integrability condition holds; local integrability of the lifted coefficients follows from smoothness. The superposition theorem~\cite{BogachevRocknerShaposhnikov2019}*{Theorem 1.1} therefore represents the stationary lifted path by a martingale solution. As in \cref{lem:superposition}, the paths remain on the closed graph of $\iota$. Projection gives a solution of the Langevin martingale problem with initial law $\mu$. Local uniqueness and nonexplosion identify it with the Langevin process, proving invariance.

Invariance of $\mu$ ensures that the density semigroup $\Pt_t$ is Markov. Path continuity and invariance imply strong continuity of the observable semigroup in $L^2(\mu)$, first on bounded continuous functions and then, by contraction and density, on all of $L^2(\mu)$. Its adjoint $\Pt_t$ is therefore strongly continuous in $L^2(\mu)$, and contraction and density extend this property to $L^1(\mu)$.

The velocity noise vector fields and their Lie brackets with the kinetic drift span the tangent space at every point. H\"ormander's local hypoellipticity theorem~\cite{Hormander1967}*{Theorem 1.1} therefore yields smooth densities at positive times. 

We prove the entropy balance using the truncations $\beta_n$ from \eqref{eq:beta}. Set $g_t=\Pt_tg$ and $\psi_R=\chi(E/R)$, with the nonincreasing cutoff from the proof of \cref{prop:entropy-chain}. Smoothness of $\Pt_tg$ and compact support of $\psi_R$ justify the chain rule and integration by parts at positive times, giving
\[
 \begin{aligned}
 \frac{\ud}{\ud t}\int\psi_R\beta_n(g_t)\ud\mu
 &+\gamma\int\psi_R\frac{\ud^2\beta_n}{\ud r^2}(g_t)|\nabla_vg_t|^2\ud\mu\\
 &=\int\beta_n(g_t)(\La+\gamma\Ls)\psi_R\ud\mu.
 \end{aligned}
\]
Both generators act on $\psi_R$ in this identity, so the cutoff error involves $\beta_n(g_t)$ but not $\nabla_vg_t$. Since $\Ls E=d-|v|^2$ and $|\nabla_vE|^2=|v|^2$, the bounds \eqref{eq:energyderivatives} give
\[
 |(\La+\gamma\Ls)\psi_R|
       \le C\mathbf1_{\{E\ge R\}},\qquad R\ge1.
\]
Integrate between $0<s<t\le T$. Since $\beta_n(g)\le C_n(g+1)$, the time integral of the cutoff error tends to zero as $R\to\infty$. The traces converge by dominated convergence, while $\psi_R\uparrow1$ permits monotone convergence in the nonnegative dissipation term. We obtain
\[
 \int\beta_n(g_t)\ud\mu
 +\gamma\int_s^t\!\int\frac{\ud^2\beta_n}{\ud r^2}(g_\tau)|\nabla_vg_\tau|^2\ud\mu\ud\tau
 =\int\beta_n(g_s)\ud\mu.
\]
Strong $L^1(\mu)$ continuity and the Lipschitz bound on $\beta_n$ allow $s\downarrow0$. Finally, $\beta_n(r)\uparrow r\log r-r+1$ and $r\frac{\ud^2\beta_n}{\ud r^2}(r)=a_n(r)\uparrow1$ for $r>0$. Monotone convergence therefore gives \eqref{eq:semigroup-entropy-balance}, with the velocity Fisher information defined by \eqref{eq:fisher-definition}. In particular, the velocity current $W_t=-\gamma\nabla_vg_t$ of the Langevin evolution has finite quadratic action on every bounded interval.

Dynkin's formula for test functions in $C_c^\infty(\RR^{2d})$, together with smoothness of $g_t$ at positive times, shows that $(g,W)$ satisfies \eqref{eq:controlled-intro} in distributions on $(0,T)\times\RR^{2d}$. The path $\nu_t=g_t\mu$ has a density at every time and is narrowly continuous on $[0,T]$ by strong $L^1(\mu)$ continuity. Since $\nu_0\in\PPtwo$ and $\Ent_\mu(g)<\infty$, \cref{lem:superposition} with $A_t=-\gamma\nabla_v\log g_t$ gives $\nu\in C([0,T];\PPtwo)$.
\end{proof}

\subsection{Approximation on bounded sets}\label{app:compact-approximation}

We now construct smooth, compactly supported approximations of the density-current pair with convergence of the transport distances, the action, and the total and spatial marginal entropies. These convergences allow us to pass to the limit in \eqref{eq:weak-comparison}.

\begin{theorem}[Approximation on bounded sets]\label{thm:recovery}
Under the assumptions of \cref{thm:weak}, there exist smooth density-current pairs $(\varrho^k,F^k)$ in Lebesgue variables satisfying \eqref{eq:ordinary-control}, with $\varrho^k\ge0$ and $\int\varrho_t^k\ud z=1$ at every time. For each $k$, both $\varrho^k$ and $F^k$ are compactly supported in phase space, uniformly in time on a neighborhood of $[0,T]$. As $k\to\infty$,
\begin{equation}\label{eq:compact-recovery}
 \begin{aligned}
 &\sup_{t\in[0,T]}W_2(\varrho_t^k\ud z,\nu_t)\longrightarrow0, \qquad \int_0^T\!\int(1+|z|^2)|\varrho^k-\varrho|\ud z\ud t
   \longrightarrow0,\\
 &\int_0^T\!\int\frac{|F^k|^2}{\varrho^k}\ud z\ud t
   \longrightarrow
   \int_0^T\!\int\frac{|F|^2}{\varrho}\ud z\ud t, \qquad \frac{F^k}{\sqrt{\varrho^k}}
   \longrightarrow\frac{F}{\sqrt{\varrho}}
   \qquad\text{in }L^2(\ud t\ud z).
 \end{aligned}
\end{equation}
The total and spatial marginal entropies also converge in $L^1(0,T)$.
\end{theorem}

\begin{proof}
Let $\mathbf P$ be the superposition measure from \cref{lem:superposition}. To obtain a compactly supported path while preserving its controlled equation, define
\[
 \Omega_R=\left\{\sup_{t\in[0,T]}|Z_t|\le R\right\},\qquad
 p_R=\mathbf P(\Omega_R)\uparrow1.
\]
For $p_R>0$, condition $\mathbf P$ on $\Omega_R$. Its marginals $\nu_t^R=\varrho_t^R\ud z$ satisfy
\[
 0\le\varrho_t^R\le p_R^{-1}\varrho_t,\qquad
 \operatorname{supp}\nu_t^R\subset\overline B_R,
 \qquad F_t^R=\varrho_t^RA_t.
\]
The conditioned trajectories still satisfy \eqref{eq:characteristics}. Testing that equation and taking conditional expectation shows that $(\varrho^R,F^R)$ solves \eqref{eq:ordinary-control}. The marginals are $W_2$-continuous, since the trajectories are continuous and remain in $\overline B_R$.

The moment bound \eqref{eq:supmoment} and finite quadratic action give, as $R\to\infty$,
\begin{equation}\label{eq:restrictionconvergence}
 \begin{aligned}
 &\sup_{t\in[0,T]}\int(1+|z|^2)|\varrho_t^R-\varrho_t|\ud z
   \longrightarrow0,\\
 &\int_0^T\!\int|A|^2\varrho^R\ud z\ud t
   \longrightarrow\int_0^T\!\int|A|^2\varrho\ud z\ud t.
 \end{aligned}
\end{equation}
Before normalization, the contribution of the discarded trajectories to the weighted density error is bounded uniformly in time by
$\mathbb E_{\mathbf P}[(1+\sup_{t\in[0,T]}|Z_t|^2)\mathbf1_{\Omega_R^c}]$.
This tends to zero by dominated convergence, while the normalization error vanishes because $p_R\to1$.
For the action convergence, finite quadratic action makes
$\int_0^T|A_t(Z_t)|^2\ud t$ integrable on path space, so monotone convergence applies to the restrictions to $\Omega_R$ before normalization.
Likewise, integrability of $\int_0^T|A_t(Z_t)|\ud t$, which follows by Cauchy--Schwarz, gives $F^R\to F$ in space-time $L^1$.

The entropy chain rule gives
$\sup_{t\in[0,T]}\Ent_\mu(g_t)<\infty$.
Set $g_t^R=\varrho_t^R/\varrho_\infty$.
Since $g_t^R\le g_t/p_R$,
\[
 (g_t^R\log g_t^R)^+
 \le p_R^{-1}g_t\bigl(\log_+g_t+\log p_R^{-1}\bigr),
 \qquad
 (g_t^R\log g_t^R)^-\le \e^{-1}.
\]
Together with density convergence and the uniform entropy bound, these estimates imply convergence of the total entropy at every fixed time and in $L^1(0,T)$.
The same argument applies to the spatial marginal entropies, since $q_t^R\le q_t/p_R$.

Fix $R$ and apply the smoothing \eqref{eq:shear-recovery} to $(\varrho^R,F^R)$, obtaining $(\varrho^{R,\varepsilon},F^{R,\varepsilon})$.
For $0<\varepsilon\le1$, the supports lie in a fixed compact set on which $\nabla U$ is Lipschitz.
Choosing $0\le\psi\le1$ equal to one on a neighborhood of this set, \eqref{eq:jensen-smoothing} gives
\[
 \left(\int\frac{|F_t^{R,\varepsilon}|^2}{\varrho_t^{R,\varepsilon}}\ud z\right)^{1/2}
 \le
 \left(\int\frac{|F_t^R|^2}{\varrho_t^R}\ud z\right)^{1/2}
 +C_R\varepsilon.
\]
Taking the $L^2(0,T)$ norm gives the upper bound for the integrated action as $\varepsilon\downarrow0$.
Distributional convergence and lower semicontinuity give the corresponding lower bound, so the integrated action converges.
The currents also converge in space-time $L^1$, while convolution gives uniform convergence in time of the laws in $W_2$ and weighted space-time $L^1$ convergence of the densities.

For the entropy, Jensen's inequality shows that convolution does not increase $\int\varrho\log\varrho\ud z$.
Since the kernel is centered, the linear term in the Taylor expansion of $U(x)+|v|^2/2$ averages to zero.
The local Hessian bound therefore yields
\begin{equation}\label{eq:compactentropy}
 \Ent(\varrho_t^{R,\varepsilon}\ud z\mid\mu)
 \le \Ent(\nu_t^R\mid\mu)+C_R\varepsilon^2.
\end{equation}
The same argument gives the analogous estimate for the spatial marginal entropy.

It remains to smooth in time while preserving the endpoint limits. For fixed $R,\varepsilon$, narrow continuity of $\nu^R$ and smoothness of the compact kernel imply that the smoothed densities and all their phase-space derivatives are continuous in time, uniformly on their compact support. Extend the endpoint densities to a slightly larger time interval by the Hamiltonian flow, and set the velocity current to zero outside $[0,T]$. By \eqref{eq:energyderivatives}, the flow is complete in both time directions and maps compact sets to compact sets over bounded intervals. The extended density is continuous across the endpoints, so the extended pair solves \eqref{eq:ordinary-control} without additional sources there.

Choose $\zeta\in C_c^\infty((-1,1))$ with $\zeta\ge0$ and $\int\zeta=1$, and set $\zeta_\delta(t)=\delta^{-1}\zeta(t/\delta)$.
For sufficiently small $\delta>0$, we consider the convolution of the extended pair in time: 
\[
 (\varrho^{R,\varepsilon,\delta},F^{R,\varepsilon,\delta})
 =\zeta_\delta*_t(\varrho^{R,\varepsilon},F^{R,\varepsilon}).
\]
Since the coefficients of \eqref{eq:ordinary-control} are independent of time, this convolution preserves the equation.
It also preserves nonnegativity and unit mass, and yields a pair that is smooth in all variables on a neighborhood of $[0,T]$.

Jensen's inequality and Fubini's theorem give
\[
 \begin{aligned}
 \int_0^T\!\int
 \frac{|F^{R,\varepsilon,\delta}|^2}
      {\varrho^{R,\varepsilon,\delta}}\ud z\ud t
 &\le
 \int_{-\delta}^{T+\delta}
 \left(\int_0^T\zeta_\delta(t-s)\ud t\right)
 \int\frac{|F_s^{R,\varepsilon}|^2}
          {\varrho_s^{R,\varepsilon}}\ud z\ud s\\
 &\le
 \int_0^T\!\int
 \frac{|F^{R,\varepsilon}|^2}
      {\varrho^{R,\varepsilon}}\ud z\ud t.
 \end{aligned}
\]
The last inequality follows from
$0\le\int_0^T\zeta_\delta(t-s)\ud t\le1$
and $F^{R,\varepsilon}=0$ outside $[0,T]$.
The approximate-identity property gives convergence of the density and current in space-time $L^1$ as $\delta\downarrow0$.
Lower semicontinuity, together with the upper bound above, therefore gives convergence of the integrated action.

Continuity of the extended density on its common compact support implies uniform convergence in $(t,z)$ as $\delta\downarrow0$.
Consequently, the densities converge in weighted $L^1$ and the corresponding laws converge in $W_2$, both uniformly in time.
Since $r\log r$, with $0\log0=0$, is uniformly continuous on bounded intervals, the total and spatial marginal entropies also converge uniformly in time.

We now combine the three approximations into a single sequence.
Let $R_k\uparrow\infty$ with $p_{R_k}>0$.
For each $k$, choose $0<\varepsilon_k\le1/k$ so that phase-space smoothing changes the integrated action, the current in space-time $L^1$, and the density in weighted space-time $L^1$ by at most $1/k$ each, and so that $C_{R_k}\varepsilon_k^2\le1/k$ in \eqref{eq:compactentropy}.
Then choose $0<\delta_k\le1/k$ so that time smoothing changes the integrated action and the current in space-time $L^1$ by at most $1/k$ each.
We also require that it changes the density in weighted $L^1$, the law in $W_2$, and both entropies by at most $1/k$ each, uniformly for $t\in[0,T]$.
Set
\[
 (\varrho^k,F^k)
 =(\varrho^{R_k,\varepsilon_k,\delta_k},
   F^{R_k,\varepsilon_k,\delta_k}),
 \qquad
 g^k=\varrho^k/\varrho_\infty.
\]

The restriction estimate \eqref{eq:restrictionconvergence}, the shrinking support of the phase-space convolution kernel, and the chosen error from time smoothing give uniform convergence in time in $W_2$.
The upper bound \eqref{eq:compactentropy}, convergence of the restricted entropies, and lower semicontinuity imply convergence of the total and spatial marginal entropies at every time.
Moreover, the uniform entropy bound for $g_t$, the restriction bound $g_t^R\le g_t/p_R$, and the chosen smoothing errors bound the approximating entropies uniformly in $k$ and $t$.
Dominated convergence therefore yields entropy convergence in $L^1(0,T)$.

Finally, the integrated action converges by construction, and both $F^k\to F$ and $\varrho^k\to\varrho$ hold in $L^1$.
Applying \cref{lem:weighted-fields} to $(\varrho^k,F^k)$ gives the strong $L^2(\ud t\ud z)$ convergence of $F^k/\sqrt{\varrho^k}$ in \eqref{eq:compact-recovery}.
\end{proof}

\subsection{Transport stability}\label{app:transport-stability}

The following lemma gives continuity of the Wasserstein current and, when the densities also converge in $L^1$, the strong convergence of the weighted Brenier displacements used to pass the control pairing to the limit. Both follow from stability of the joint law of position, velocity, and transported position. The optimal transportation maps depend only on the spatial laws, so the lemma applies with any reference measure used to express the densities.

\begin{lemma}[Transport stability]\label{lem:transport-stability}
Let $\sigma\in\PPtwo(\RR^d)$. Let $\nu_n,\nu\in\PPtwo(\RR^{2d})$ have absolutely continuous spatial marginals $q_n\mux,q\mux$, and suppose $\nu_n\to\nu$ in $W_2$. Then the joint laws
\[
 (x,v,T_{q_n,\sigma}(x))_\#\nu_n
\]
converge in $W_2$ to $(x,v,T_{q,\sigma}(x))_\#\nu$. In particular, $\Cs(\nu_n)\to\Cs(\nu)$.

Let now $\nu^n,\nu\in C([0,T];\PPtwo(\RR^{2d}))$ have Lebesgue densities $\varrho^n,\varrho$, and suppose $\sup_{t\in[0,T]}W_2(\nu_t^n,\nu_t)\to0$ and $\varrho^n\to\varrho$ in $L^1((0,T)\times\RR^{2d})$. Then
\begin{equation}\label{eq:displacement-strong}
 \sqrt{\varrho^n}\,\xi_{q^n,\sigma}
 \longrightarrow\sqrt\varrho\,\xi_{q,\sigma}
 \quad\text{strongly in }L^2(\ud t\ud z).
\end{equation}
\end{lemma}
\begin{proof}
The joint laws are tight because their coordinate marginals converge. Stability of optimal plans~\cite{ambrosio2005gradient}*{Proposition 7.1.3} identifies the $(x,y)$ marginal of any weak limit with the unique optimal plan from $q\mux$ to $\sigma$. This plan is concentrated on the graph $y=T_{q,\sigma}(x)$. Since the $(x,v)$ marginal is $\nu$, the full joint limit is uniquely determined. The second moments converge as the sum of the second moments of the three coordinate marginals. Hence the joint laws converge in $W_2$.

The Wasserstein current is the integral of $v\cdot(x-y)$ against the joint law. This continuous function has at most quadratic growth, so joint $W_2$ convergence gives convergence of the currents.

For \eqref{eq:displacement-strong}, convergence of the joint laws at every time identifies the distributional limit of $\varrho^n\xi_{q^n,\sigma}$ on space-time as $\varrho\xi_{q,\sigma}$. The required domination in time follows from the uniform second-moment bounds. The triangle inequality gives uniform convergence of $W_2(q_t^n\mux,\sigma)$ to $W_2(q_t\mux,\sigma)$, so
\[
\int_0^T\!\int\varrho^n|\xi_{q^n,\sigma}|^2\ud z\ud t =\int_0^T W_2^2(q_t^n\mux,\sigma)\ud t \longrightarrow\int_0^T W_2^2(q_t\mux,\sigma)\ud t.
\]
Apply \cref{lem:weighted-fields} with $r_n=\varrho^n$ and $B_n=\varrho^n\xi_{q^n,\sigma}$. The assumed $L^1$ convergence of the densities and the convergence of the quadratic integrals prove \eqref{eq:displacement-strong}.
\end{proof}

\end{document}